\documentclass[12pt,reqno]{amsart}
\usepackage{mathrsfs}
\usepackage{amsmath}
\usepackage{amssymb, amsmath}

\newcommand{\newcom}{\newcommand}

\newcom{\al}{\alpha}
\newcom{\be}{\beta}
\newcom{\eps}{\epsilon}
\newcom{\ga}{\gamma}
\newcom{\Ga}{\Gamma}
\newcom{\ka}{\kappa}
\newcom{\Lam}{\Lambda}
\newcom{\lam}{\lambda}
\newcom{\Om}{\Omega}
\newcom{\om}{\omega}
\newcom{\Si}{\Sigma}
\newcom{\si}{\sigma}
\newcom{\tht}{\theta}
\newcom{\dtri}{\nabla}
\newcom{\tri}{\triangle}
\newcom{\oo}{\infty}
\newcom{\vphi}{\varphi}
\newcom{\calb}{{\mathcal B}}
\newcom{\calc}{{\mathcal C}}
\newcom{\cD}{{\mathcal D}}
\newcom{\cF}{{\mathcal F}}
\newcom{\cI}{{\mathcal I}}
\newcom{\cL}{{\mathcal L}}
\newcom{\cM}{{\mathcal M}}
\newcom{\cP}{{\mathcal P}}
\newcom{\cR}{{\mathcal R}}
\newcom{\cS}{{\mathcal S}}
\newcom{\cQ}{{\mathcal Q}}
\newcom{\caly}{{\mathcal Y}}
\newcom{\calz}{{\mathcal Z}}
\newcom{\bfz}{{\bf Z}}
\newcom{\R}{\Bbb R}
\newcom{\N}{\Bbb N}
\newcom{\Z}{\Bbb Z}
\newcom{\C}{\Bbb C}
\newcom{\E}{\Bbb E}

\newcom{\f}{\frac}
\newcom{\di}{\displaystyle\int}
\newcom{\ds}{\displaystyle\sum}
\newcom{\dl}{\displaystyle\lim}
\newcom{\ov}{\overline}
\newcom{\sset}{\subset}
\newcom{\wt}{\widetilde}
\newcom{\p}{\partial}
\newcom\na{\nabla}
\newcom{\co}{\cdot}
\newcom{\suml}{\sum\limits}
\newcom{\supl}{\sup\limits}
\newcom{\intl}{\int\limits}
\newcom{\infl}{\inf\limits}
\newcom{\disp}{\displaystyle}
\newcom{\non}{\nonumber}
\newcom{\no}{\noindent}
\newcom{\QED}{$\square$}
\def\ef{\hphantom{MM}\hfill\llap{$\square$}\goodbreak}

\newtheorem{athm}{\bf \t}[section]
\newenvironment{thm} [1] {\def\t{#1}\begin{athm} \bf \rm} {\end {athm}}
\newcom{\bthm}{\begin{thm}}\newcom{\ethm}{\end{thm}}

\newcom{\beq}{\begin{equation}}
\newcom{\eeq}{\end{equation}}
\newcom{\ben}{\begin{eqnarray}}
\newcom{\een}{\end{eqnarray}}
\newcom{\beno}{\begin{eqnarray*}}
\newcom{\eeno}{\end{eqnarray*}}

\numberwithin{equation}{section}

\begin{document}

\title[water-wave problem with surface tension]
{well-posedness of the water-wave problem with surface tension}

\author{Mei Ming}
\address{Academy of
Mathematics $\&$ Systems Science, CAS, Beijing 100190, P. R. China}
\email{mingmeim2@gmail.com}

\author{Zhifei Zhang}
\address{School of Mathematical Sciences, Peking University, 100871, P. R. China and
Departement de Mathematique, Universite Paris Sud, 91405 Orsay Cedex, FRANCE }
\email{zfzhang@math.pku.edu.cn}

\thanks{The second author is partially supported by NSF of China under Grant 10601002}

\date{June 26}



\begin{abstract}
In this paper, we prove the local well-posedness of the water wave problem with surface tension
in the case of finite depth by working in the Eulerian setting. For the flat bottom, as surface tension tends to zero,
the solution of the water wave problem with surface tension converges to the solution of the water wave problem
without surface tension.
\end{abstract}

\maketitle

\section{Introduction}

\subsection{Presentation of the problem}

In this paper, we are concerned with the motion of an ideal,
incompressible, irrotational gravity fluid influenced by surface
tension on its surface in the case of finite depth. We restrict our
attention to the case when the surface is a graph parameterized by a
function $\zeta(t,X)$ where $t$ denotes the time variable, and
$X=(X_1,...,X_d)$ denotes the horizontal spacial variables. The
bottom of fluid is parameterized by a function $b(X)$. We denote the
fluid domain at time $t$ by $\Omega_t $. The motion of the fluid in
$\Omega_t$ is described by the incompressible Euler equation
\beq\label{1.1}
\p_t V+V\cdot \na_{X,y} V=-ge_{d+1}-\na_{X,y} P\qquad \hbox{in}\quad \Omega_t,\,\, t\geq 0
\eeq
where $-ge_{d+1}=(0,\cdots,0,-g)$
denotes the acceleration of gravity and $V=(V_1,...,V_d,V_{d+1})$
denotes the velocity field ($V_{d+1}$ is the vertical component).
The incompressibility of the fluid is expressed by
\beq \label{1.2}
\textrm{div}\, V=0 \qquad \hbox{in}\quad \Omega_t,\,\, t\geq 0,
\eeq
and the irrotationality means that
\beq \label{1.3}
\textrm{curl}\, V=0 \qquad \hbox{in}\quad \Omega_t,\,\, t\geq 0.
\eeq

Assume that no fluid particles are transported across the surface. At the bottom, this is given by
\beq \label{1.4}
V_n|_{y=b(X)}:=\textbf{n}_-\cdot\ V|_{y=b(X)}=0 \qquad \hbox{for}\quad t>0, X\in\R^d
\eeq
where $\textbf{n}_-:=\frac{1}{\sqrt{1+|\na_X b|^2}}(\na_X b,-1)^T$ denotes the outward normal vector to
the lower boundary of $\Omega_t $. At the free surface, the boundary condition is kinematic and is given by
\beq\label{1.5}
\p_t \zeta-\sqrt{1+|\na_X \zeta|^2}V_n|_{y=\zeta(X)}=0
\qquad \hbox{for}\quad t>0, X\in\R^d\eeq
where $V_n=\textbf{n}_+\cdot V|_{y=\zeta(X)}$, with $\textbf{n}_+:=\frac{1}{\sqrt{1+|\na_X
\zeta|^2}}(-\na_X \zeta, 1)^T$ denoting the outward normal vector to the free surface.

With surface tension, the pressure at the surface is given by
\beq\label{1.6} P|_{y=\zeta(t,X)}=-\kappa\na_X\cdot (\frac{\na_X
\zeta}{\sqrt{1+|\na_X \zeta|^2}})\qquad \hbox{for}\quad t\geq 0,
X\in \R^d, \eeq
where $\kappa>0$ is the surface tension coefficient.

The above problem is known as the water wave problem. Concerning 2-D water wave problem, when
surface tension is neglected and the motion of free surface is a
small perturbation of still water, one could check Nalimov
\cite{nal}, Yosihara \cite{yos} and W. Craig \cite{wal}. In general,
the local well-posedness of the water wave problem without surface tension
was solved by S. Wu \cite{Wu1,Wu2} in the case of infinite depth and see also Ambrose and Masmoudi
\cite{amb1,amb2} for a different proof. More recently, D. Lannes \cite{lan1} considered the water-wave problem
without surface tension in the case of finite depth by working in the Eulerian setting.
We should mention some recent results concerning the rotational water wave problem\cite{lin,Co-Sh,Sh-Ze, zzc}.

The purpose of this paper is to study the water wave problem with surface tension
and zero surface tension limit in the case of finite depth.
Although it seems possible to adapt the method of Ambrose and Masmoudi\cite{amb1,amb2} to the case of finite depth,
we choose to work in the Eulerian setting as in \cite{lan1}, since it's the easiest to handle,
especially when the asymptotic properties of the solutions are concerned\cite{Alv1,Bon}.
On the other hand, the water wave problem can be reformulated as a Hamiltonian system
in the Eulerian coordinates\cite{zak}. To exploit the theory(for example, canonical transformation theory)
in the Hamiltonian mechanics, it is natural to use the coordinates in which the water wave problem
has a Hamiltonian structure.

\subsection{Presentation of the result}

As in \cite{lan1}, we use an alternative formulation of the water wave system (\ref{1.1})-(\ref{1.6}).
From (\ref{1.2}) and (\ref{1.3}), there exsits a
potential flow function $\phi$ such that $V=\na_{X,y} \phi$ and
\beq\label{1.7}
\Delta_{X,y}\phi=0 \qquad \hbox{in}\quad \Omega_t,\,\, t\geq 0.
\eeq
The boundary conditions (\ref{1.4}) and (\ref{1.5}) can be expressed in terms of $\phi$
\beq \label{1.8}
\p_{n_-} \phi|_{y=b(X)}=0, \qquad \hbox{for}\quad t>0,\,\, X\in\R^d,
\eeq
and
\beq \label{1.9}
\p_t \zeta-\sqrt{1+|\na_X \zeta
|^2}\p_{n_+} \phi|_{y=\zeta(X)}=0, \qquad \hbox{for}\quad t>0,
\,\,X\in\R^d, \eeq
where we denote $\p_{n_-}:=\textbf{n}_-\cdot \na_{X,y}$ and
$\p_{n_+}:=\textbf{n}_+\cdot \na_{X,y}$. The Euler's
equation (\ref{1.1}) can be put into Bernoulli's form
\beq\label{1.10}
\p_t \phi+\frac{1}{2}|\na_{X,y} \phi|^2+gy=-P \qquad \hbox{in}\quad \Omega_t,\,\, t\geq 0.
\eeq

We next reduce the system (\ref{1.7})-(\ref{1.10}) to a system where all the functions are evaluated at the free surface only.
For this purpose, we introduce the trace of the velocity potential $\phi$
at the free surface
\[ \psi(t,X):=\phi(t,X,\zeta(t,X)), \]
and the (rescaled) Dirichlet-Neumann operator $G(\zeta,b)$(or simply
$G(\zeta)$)
\beno
G(\zeta)\psi:=\sqrt{1+|\na_X\zeta|^2}\p_{n_+}\phi|_{y=\zeta(t,X)}.
\eeno

Taking the trace of (\ref{1.10}) on the free surface, the
system (\ref{1.7})-(\ref{1.10}) is equivalent to the system
\beq \label{1.12}
\left\{
\begin{array}{ll}
\p_t \zeta-G(\zeta)\psi=0, \\
\p_t \psi+g\zeta+\frac{1}{2}|\na_X
\psi|^2-\frac{(G(\zeta)\psi+\na_X\zeta\cdot \na_X
\psi)^2}{2(1+|\na_X \zeta|^2)}=\kappa\na_X \cdot (\frac{\na_X
\zeta}{\sqrt{1+|\na_X \zeta|^2}}),
\end{array}\right.
\eeq
which is an evolution equation for the height of the free
surface $\zeta(t,X)$ and the trace of the velocity potential on the
free surface $\psi(t,X)$. Our results in this paper will be given for this system.

We will essentially follow the framework of \cite{lan1}. It is well-known that Taylor's
sign condition is a necessary condition for well-posedness of the water-wave system
without surface tension\cite{Ebin,htz}. It is unnecessary for the system with surface tension
due to the ``smoothing effect" of surface tension. Due to the derivative loss in the system (\ref{1.12}),
we will use Nash-Moser iteration to solve (\ref{1.12}).
The key step is to obtain the well-posedness
and the tame estimates of the soluion to the linearized equations of (\ref{1.12}):
\ben\label{linear} \ \left\{
\begin{array}{ll}
\p_tV_1+\na_X\cdot(\underline{v}V_1)-\underline{G}V_2=H_1, \\
\p_tV_2+(\underline{a}-\underline{A})V_1+\underline{v}\cdot \na_X V_2=H_2,
\end{array}\right.
\een
where $\underline{a}, \underline{v}$ are smooth functions, the source terms $(H_1,H_2)\in H^{k+1}\times H^{k+\f12},
\underline{G}=G(\underline{\zeta})$ and
\[
\underline{A}=\kappa\na_X\cdot\Bigl[ \frac{\na_X}{\sqrt{1+|\na_X
\underline{\zeta}|^2}}-\frac{\na_X
\underline{\zeta}(\na_X\underline{\zeta}\cdot \na_X)}{(1+|\na_X
\underline{\zeta}|^2)^\frac{3}{2}}\Bigr]. \]
As usual(find a change of unknowns which symmetrizes  the system (\ref{linear})), if we try to use the energy functional
$$
E_k(V)=(\Lambda^kV_1, -\underline{A}\Lambda^kV_1)+(\Lambda^kV_2, \underline{G}\Lambda^kV_2), \quad \Lambda=(1+|D|^2)^\f12
$$
we have to deal with some singular terms like
$$
(\Lambda^kV_2, [\underline{G},-\underline{A}]\Lambda^kV_1).
$$
Note that $V_1\in H^{k+1}$ and $V_2\in H^{k+\f12}$.
It seems impossible to control this term by $|V_1|_{H^{k+1}}|V_2|_{H^{k+\f12}}$,
since the commutator $[\underline{G},-\underline{A}]$ is an operator of order $2$.
Fortunately, we find an important fact that the main part of the Dirichlet-Neumann operator $\underline{G}$ is similar to
that of the operator $\underline{A}$. Formally, if the main part of $\underline{G}$ is the same as
that of the operator $\underline{A}$, we expect that $[\underline{G},-\underline{A}]$ becomes an operator of order $1$.
In order to use this fact, we need to introduce a more complicated energy functional
\[E_k(V):=(\tilde{\Lambda}^k\, \sigma V_1
,\sigma^{-1}(\underline{a}-\underline{A})\sigma^{-1}\tilde{\Lambda}^k\,\sigma
V_1)+(\tilde{\Lambda}^k \,\sigma^{-1}V_2,\,\,\sigma \underline{G}
\,\sigma \tilde{\Lambda}^k \,\sigma^{-1} V_2),
\]
where
\beno
&&\widetilde{\Lambda}=|D|^2-\frac{\p_i\underline{\zeta}\p_j\underline{\zeta}}
{1+|\na_X \underline{\zeta}|^2}D_iD_j,\quad\sigma=(1+|\na_X
\underline{\zeta}|^2)^{-\frac{1}{4}}.
\eeno
On the other hand, to obtain the tame estimates of the solution,
we rely heavily on some sharp pseudo-differential operator estimates obtained in \cite{lan2}.

Now we state our results as follows.

\bthm{Theorem}\label{thm1.1}
Let $b\in C^\infty_b(\R^d)$. There exists $P>D>0$ such that for all $\zeta_0\in H^{s+P}(\R^d)$ and $\psi_0$
be such that  $\na_X \psi_0\in H^{s+P}(\R^d)^d$, with $s>M$ ($M$ depends
on $d$ only). Assume moreover that
\[
\min\{\zeta_0-b,-b\}\ge 2h_0 \qquad\hbox{on}\quad \R^d \qquad\hbox{for
some} \quad h_0\ge 0.
\]
Then there exist $T>0$ and a unique solution $(\zeta,\psi)$ to
the water-wave system (\ref{1.12}) with the initial condition
$(\zeta_0,\psi_0)$ and such that $(\zeta,\psi-\psi_0)\in
C^1([0,T],\,H^{s+D}(\R^d)\times H^{s+D}(\R^d))$.
\ethm

For the flat bottom, we also obtain the zero surface tension limit.

\bthm{Theorem}\label{thm1.2} Let $b=-1$ and $(\zeta_0, \psi_0)$ satisfy the same assumptions as Theorem \ref{thm1.1}.
Assume that the surface tension coefficient $\kappa$ is sufficiently small.
Then there exist $T>0$ independent of $\kappa$  and a unique solution $(\zeta^\kappa,\psi^\kappa)$ to
(\ref{1.12}) with $\kappa>0$  such that $(\zeta^\kappa,\psi^\kappa-\psi_0)\in
C^1([0,T],\,H^{s+D}(\R^d)\times H^{s+D}(\R^d))$.
Moreover, as $\kappa$ tends to zero, the solution $(\zeta^\kappa,\psi^\kappa)$ converges to the solution
$(\zeta,\psi)$ of (\ref{1.12}) with $\kappa=0$.
\ethm

\bthm{Remark} In \cite{lan1}, D. Lannes proved the well-posedness of the water wave problem without surface tension
under the following condition on the bottom
\ben\label{bottom}
\Pi_b(V_0|_{y=b(x)},V_0|_{y=b(x)})\le \f g {\sqrt{1+|\na_X b|^2}},
\een
where $\Pi_b$ is the second fundamental form of $\Gamma_b$, and $V_0$ is the velocity field associated to $\psi_0$.
In the case of nonzero surface tension, we don't need this extra condition (\ref{bottom}).
\ethm

\bthm{Remark} Under the condition (\ref{bottom}), D. Lannes proved the L\'{e}vy
condition in the case of zero surface tension. However, we don't know whether the L\'{e}vy
condition still holds for small surface tension under (\ref{bottom}). This is the reason why we
restrict Theorem \ref{thm1.2} to the case of the flat bottom.
\ethm

\noindent{\bf Organization of the paper.} In section 2, we prsent a tame elliptic estimate on a strip.
In section 3, we review some sharp pseudo-differential operator estimates,
and then introduce an important elliptic operator and study its properties.
Section 4 is devoted to the study of the Dirichlet-Neumann operator.
In section 5, we prove the well-posedness of the
linearized water-wave equations and the tame energy estimates of the solution. In section 6, we
solve the fully nonlinear equations by using  Nash-Moser iteration and study zero surface tension limit.

\subsection{Notation}

We list some notations we will use throughout this paper:

-$C$ denotes some numerical constant which  may change from one line  to another. If
the constant $C$ depends on $\lambda_1,\lambda_2,...$, we simply denote it by $C(\lambda_1,\lambda_2,...)$.

-We denote by $m_0$ the first integer strictly larger then $\f {d+1} 2$.

-We write $|\alpha|=\alpha_1+...+\alpha_{d+1}$ with
$\alpha=(\alpha_1,...,\alpha_{d+1})$.

-We denote $\p_i=\p_{X_i}$, for $i=1,...,d$, $\p_{d+1}=\p_y$, and
for $\alpha\in \N^{d+1},
\p^{\alpha}={\p_1}^{\alpha_1}...{\p_{d+1}}^{\alpha_{d+1}}$.

-We denote by $C^k_b(\R^d) $ the set of bounded and continuous on $\R^d$ functions together
with their derivatives of order less than or equal to $k$. We also
denote $C^\infty_b=\cap_k C^k_b$.

-We denote by $(\cdot,\cdot)$ the usual real scalar product on
$L^2(\R^d)$.

-We denote by $\Lambda=\Lambda(D)$ the Fourier multiplier with the
symbol $\Lambda(\xi)=(1+|\xi|^2)^{\frac{1}{2}}$.

-We denote by $H^s(\R^d)$ the Sobolev space with norm
$|f|_{H^s}:=(\int \Lambda(\xi)^{2s}|\hat{f}(\xi)|^2
d\xi)^{\frac{1}{2}}$.

-We denote $|f|_{H^s_T}=\sup_{t\in[0,T]}|f(t)|_{H^s}$, when $f\in
C([0,T],H^s)$.

-We denote $|F|_B=|f_1|_B+...+|f_n|_B$, where $F=(f_1,...,f_n)$, $B$
is a Banach space, and $F\in B^n$.

-For an open set $U\subset \R^{d+1}$ we denote by $\|\cdot:U\|_p$,
$\|\cdot:U\|_{k,\infty}$, and $\|\cdot:U\|_{k,2}$ the norms of
$L^p(U)$, $W^{k,\infty}(U)$ and $H^k(U)$ respectively. When no
confusion can be made we omit $U$.

\section{Elliptic boundary value problem on a strip}
Assume that  $\Omega=\{(X,y)\in \R^{d+1}, b(X)<y<a(X)\}$, where
$a(X)$ and $b(X)$ satisfy \beq\label{2.1} \exists h_0>0, \quad
\hbox{s.t. on}\,\,\R^d,\,\,\min\{-b(X),a(X)-b(X)\}\ge h_0. \eeq
We study the boundary value problem of the form \beno \ \left\{
\begin{array}{ll}
-\Delta_{X,y} u=h \qquad \hbox{on}\quad \Omega, \\
u|_{y=a(X)}=f,\,\,\p_{n_-}u|_{y=b(X)}=g
\end{array}\right.
\eeno
where $h$ is a function defined on $\Omega$, $f,g$ are
functions defined on $\R^d$, and $\p_{n_-}u|_{y=b(X)}$ denotes
the normal derivative of $u$ at the boundary $y=b(X)$.

We denote by $\cS=\R^d\times (-1,0)$ a flat strip. We denote by $S$ a
diffeomorphism from $\cS$ to $\Omega$, which is of the form
\[
S:\begin{array}{l}\cS\longrightarrow \Omega\\
(\tilde{X},\tilde{y})\mapsto
(\tilde{X},s(\tilde{X},\tilde{y})),\end{array}
\]
where $s(\widetilde{X},\widetilde{y})=-b(\widetilde{X})\widetilde{y}+(1+\widetilde{y})a(\widetilde{X}).
$

Using the diffeomorphism $S$, the elliptic equation $-\Delta_{X,y} u=0$ on $\Omega$
can be equivalently formulated as a variable coefficients equation $\widetilde{\textbf{P}}\widetilde{u}=0$ on $\cS$.
More precisely,

\bthm{Lemma}\label{flat-ellip} The elliptic equation $-\Delta_{X,y} u=h$ holds in $\cD'(\Omega)$
if and only if  the equation $\tilde{\textbf{P}}\tilde{u}=(a-b)\tilde{h}$
holds in $\cD'(\cS)$, where $\tilde{u}=u\circ S$ and
$\tilde{h}=h\circ S$, and
$\tilde{\textbf{P}}:=-\na_{X,y}\cdot\tilde{P}\na_{X,y}$, with
\[
\tilde{P}=\frac1{a-b}\left(\begin{matrix}
(a-b)^2I_{d\times d} & (b-a)\na_{\widetilde{X}}s\\
(b-a)\na_{\widetilde{X}}s^T & 1+|\na_{\widetilde{X}}s|^2\end{matrix}\right),\quad
\na_{\widetilde{X}}s=-\na_{\widetilde{X}}b\widetilde{y}+(1+\widetilde{y})\na_{\widetilde{X}}a.
\] Moreover, one has for all $\Theta\in\R^{d+1}$,
\[
\tilde{P}\Theta\cdot
\Theta\ge\tilde{p}|\Theta|^2,\qquad\hbox{with}\quad
\tilde{p}=C\frac{h^2_0}{\|a-b\|_\infty
(1+\|(\na_{\tilde{X}}a,\na_{\tilde{X}}b)\|^2_\infty)}.
\]
\ethm
The following tame estimate of the coefficient matrix $\widetilde{P}$ will be frequently used in the subsequence:
for any $k\in \N$
\ben\label{coeff}
\|\widetilde{P}\|_{k,2}\le C(h_0, |b|_{W^{k+1,\infty}}, |a|_{H^{m_0+1}})|a|_{H^{k+1}},
\een
which can be deduced by using H\"{o}lder inequality and interpolation argument (see also  Appendix A in \cite{lan1}).

We next present the tame elliptic estimate of the following variable coefficients elliptic equation:
\beq\label{2.3}
\left\{\begin{array}{l} \widetilde{\textbf{P}} u:=-\na_{X,y}\cdot \widetilde{P}\na_{X,y} u=h
\qquad\hbox{on}\quad \cS,\\
u|_{y=0}=f,\quad \p^{\widetilde{P}}_nu|_{y=-1}=g,
\end{array}\right.
\eeq
where $\p^{\widetilde{P}}_n$ denotes the conormal derivative associated to $\widetilde{P}$
\beno
\p^{\widetilde{P}}_n u|_{y=-1}=e_{d+1}\cdot \widetilde{P}\nabla_{X,y}u|_{y=-1}.
\eeno

\bthm{Proposition}\label{tame-elliptic}
Let $k\in \N$. Let $f\in H^{k+\frac32}(\R^d)$, $g\in
H^{k+\frac12}(\R^d)$ and $h\in H^k(\cS)$. If $b \in
W^{k+2,\infty}(\R^d)$ and $a \in H^{k+2}\cap
H^{2m_0+2}(\R^d)$, there exists a unique solution $u\in H^{k+2}(\cS)$ to
(\ref{2.3}). Moreover,
\begin{eqnarray*}
\|u\|_{k+2,2}&\le&
C_k(\|h\|_{k,2}+|f|_{H^{k+\frac32}}+|g|_{H^{k+\frac12}})\\
&& +C_k(\|h\|_{m_0-1,2}+|f|_{H^{m_0+\frac12}}+|g|_{H^{m_0-\frac12}})
|a|_{H^{k+2}},
\end{eqnarray*} where $C_k=C(\widetilde{p},
|b|_{W^{k+2,\infty}},|a|_{H^{2m_0+2}})$.
\ethm

\no{\bf Proof.} This is a direct corollary of (\ref{coeff}) and Theorem 2.9 in \cite{lan1} .\ef

\section{Sharp pseudo-differential operator estimates}

Let us firstly recall some tame pseudo-differential operator estimates from \cite{lan2} which play an important role
in the energy estimates.

\bthm{Definition}\label{def3.1}
Let $m\in \R$, $p\in \N$ and let $\Sigma$ be a function
defined over $\R^p_v\times \R^d_\xi$. We say that $\Sigma\in
C^\infty(\R^p,\mathcal{M}^{m})$ if

(1) $\Sigma|_{\R^p\times\{|\xi|\le1\}}\in C^\infty(\R^p;
L^\infty(|\xi|\le 1))$;

(2) for any $\alpha\in \N^p, \beta\in \N^d$, there exists a
nondecreasing function $C_{\alpha,\beta}(\cdot)$ such that
\[
\sup_{|\xi|\ge \frac14
}(1+|\xi|^2)^{(|\beta|-m)/2}\left|\p^\alpha_v
\p^\beta_\xi\Sigma(v,\xi)\right|\le C_{\alpha,\beta}(|v|).
\]
We say that $\Sigma\in C^\infty(\R^p,\mathcal{M}^{m})$ is $k$-regular at the origin
if $\Sigma|_{\R^p\times\{|\xi|\le1\}}\in C^\infty(\R^p;
W^{k,\infty}(\{|\xi|\le 1\}))$.
\ethm

We have the following tame pseudo-differential operator estimates:
\bthm{Proposition}\label{prop2.8}
Let $m\in \R$, $p\in \N$ and $d/2<t_0\le s_0$. Assume that
$\sigma(x,\xi)=\Sigma(v(x),\xi)$ with $\Sigma\in
C^\infty(\R^p,\mathcal{M}^{m})$ and $v\in H^{s_0}(\R^d)^p$. Then there hold
for all $s\in \R$ such that $\max\{-t_0,-t_0-m_1\}<s\le s_0+1$
\beno
&&|\sigma(x,D)u|_{H^s}<C_{\Sigma}(|v|_\infty)|v|_{H^{t_0}}|u|_{H^{s+m}},\quad \forall -t_0<s<t_0,\\
&&|\sigma(x,D)u|_{H^s}<C_{\Sigma}(|v|_\infty)(|v|_{H^s}|u|_{H^{m+t_0}}+|u|_{H^{s+m}}),\quad \forall\,t_0\le s\le s_0.
\eeno
\ethm

For $n\in \N$, we define $\sigma_1\sharp_n \sigma_2$ as
\beno
\sigma_1\sharp_n \sigma_2(x,\xi)=
\sum_{|\alpha|\le n}\f{(-i)^{|\alpha|}}{\alpha!}\partial_{\xi}^\alpha\sigma_1(x,\xi)\partial_x^\alpha\sigma_2(x,\xi)
\eeno
and the Poisson bracket $\{\sigma_1,\sigma_2\}_n$ as
\beno
\{\sigma_1,\sigma_2\}_n=\sigma_1\sharp_n \sigma_2(x,\xi)-\sigma_2\sharp_n \sigma_1(x,\xi).
\eeno

We have the following tame composition and commutator estimates:

\bthm{Proposition}\label{prop2.9}
Let $m_1$, $m_2\in \R$, $m:=m_1\wedge m_2$, and $d/2<t_0\le
s_0$. Let $\sigma^j(x,\xi)=\Sigma^j(v_j(x),\xi)$ with $p^j\in
\N$, $\Sigma^j\in C^\infty(\R^{p_j},\mathcal{M}^{m_j})$ and $v_j\in
H^{s_0+m\wedge n+1}(\R^d)^{p_j}$ ($j=1,2$). Assume moreover that
$\Sigma^1$ and $\Sigma^2$ are $n$-regular at the origin. Then
for all $s\in \R$ such that $\min\{-t_0,-t_0-m_1,-t_0-m_2\}\le
s\le s_0+1$ the following holds (writing $v:=(v_1,v_2)$)
\beno
&&|\textrm{Op}(\sigma^1)\circ \textrm{Op}(\sigma^2)u-\textrm{Op}(\sigma^1\#_n\sigma^2)u|_{H^s}\le
C(|v|_{W^{n+1,\infty}})\bigl[|u|_{H^{s+m_1+m_2-n-1}}\\
&&\qquad+(|v^1|_{H^{t_0+1}}|v^2|_{H^{s_+ +m\wedge
n}}+|v^2|_{H^{t_0+1}}|v^1|_{H^{s_+ +m\wedge n}})|u|_{H^{m+t_0}}\bigr].
\eeno
In particular, one has
\beno
&&|[\textrm{Op}(\sigma^1),\,\,\textrm{Op}(\sigma^2)]u-\textrm{Op}(\{\sigma^1,\,\sigma^2\}_n)u|_{H^s}
\le C(|v|_{W^{n+1,\infty}})\bigl[|u|_{H^{s+m_1+m_2-n-1}}\\
&&\qquad+(|v^1|_{H^{t_0+1}}|v^2|_{H^{s_+ +m\wedge
n}}+|v^2|_{H^{t_0+1}}|v^1|_{H^{s_+ +m\wedge n}})|u|_{H^{m+t_0}}\bigr].
\eeno
Here $s_+:=\max\{s,0\}$, $m\wedge n:=\max\{m,n\}$.
\ethm

\bthm{Remark} If the symbol $\sigma=\Sigma(v(X,y),\xi)$ and $u=u(X,y)$ for $(X,y)\in \cS, \xi\in \R^d$,
Prop. \ref{prop2.8} and \ref{prop2.9} with $H^k(\R^d)$ norm replaced by $H^k(\cS)$  norm
still remain true. We only need to notice the following two facts: Firstly, since $\cS$ is a flat strip,
there exists a linear extension operator
$E:H^k(\cS)\rightarrow H^k(\R^{d+1})$ such that(see \cite{adam})
\beno
&&Eu(X,y)=u(X,y), \quad\textrm{ a.e.}\,\, (X,y)\in \cS\quad \textrm{and}\\
&&\|Eu\|_{H^k(\R^{d+1})}\le C\|u\|_{H^k(\cS)}.
\eeno
Secondly,
\beno
\widetilde{\sigma}(X,y,D)(Eu)(X,y)=\sigma(X,y,D)u(X,y),\quad \textrm{a.e.}\,\, (X,y)\in \cS,
\eeno
where $\widetilde{\sigma}(X,y,\xi)=\Sigma((Ev)(X,y),\xi)$.
\ethm

We next introduce an elliptic differential operator $\Lambda_a$ defined by
\[
\Lambda_a=|D|^2-\frac{\p_ia\p_ja}
{1+|\na_X a|^2}D_iD_j,\quad D_j=\f {\partial_j} i,\, a\in H^\infty(\R^d)
\]
which will appear in the energy functional. Later, we will find that
the operator $\Lambda_a$ is similar to the main part of the Dirichlet-Neumann operator $G(a,b)$. Here
the repeated index denotes the summation.
In what follows, we denote by $C_k(s)$ a constant depending on $k$ and $|a|_{H^{s}}$.

\bthm{Lemma}\label{lem2.10} Let $k\in\N, s\ge 0, f\in H^{2k+s}\cap H^{m_0}(\R^d)$. Then there holds
\begin{eqnarray*}
|\Lambda_a^kf|_{H^s}\le
C_k(m_0+1)\bigl(|f|_{H^{2k+s}}+|f|_{H^{m_0}}|\na_X
a|_{H^{2k+s}}\bigr).
\end{eqnarray*}
\ethm

\noindent{\bf Proof.}\,Since lemma can be easily proved by using H\"{o}lder inequality
and interpolation argument(see Appendix A in \cite{lan1}), we omit its proof here.\ef

\bthm{Lemma}\label{lem2.11} Let $k\in \N, s\in [0,1], f\in H^{2k+s}\cap H^{m_0}(\R^d)$. Then there holds
\[
|\Lambda_a^kf|_{H^{s}}\ge
C_k(m_0+2)^{-1}|f|_{H^{2k+s}}-C_k(m_0+2)(1+|\na_X a|_{H^{2k+s}})|f|_{H^{m_0}}.
\]
\ethm

\noindent{\bf Proof.}\,We use an inductive argument on $k$. Let us firstly prove the case of $k=1$.
We can rewrite $\Lambda_a$ as
\begin{eqnarray*}
\Lambda_a&
=&[\delta_{ij}-\,(1+|\na_X
a|^2)^{-1}\p_i a\p_j
a]D_iD_j\\
&\triangleq&g_{ij}(\na_X a)D_iD_j.
\end{eqnarray*}
Then we have
\ben
(\Lambda_a f,f)&=&-(g_{ij}(\na_Xa)\p_i\p_jf,f)\nonumber\\
&=&(g_{ij}(\na_Xa)\p_if,\p_jf)+(\p_jg_{ij}(\na_Xa)\p_if,f)\nonumber\\
&\ge& C(m_0+1)^{-1}|\na f|_{L^2}^2-C(m_0+2)|f|_{L^2}|\na f|_{L^2}\nonumber\\
&\ge& C(m_0+1)^{-1}|f|_{H^1}^2-C(m_0+2)|f|_{L^2}^2.\label{2.8}
\een
Note that
\begin{eqnarray*}
|(\Lambda_a\na_X f,\na_X f)|& \le&|(\na_X\Lambda_af,\na_X f)|
+|((\na_X g_{ij})D_iD_jf,\na_X f)|\\
&\le&
(|\Lambda_af|_{L^2}+C(m_0+2)|f|_{H^1})|\na_X f|_{H^1},
\end{eqnarray*}
which together with (\ref{2.8}) gives
\begin{eqnarray*}
&&C(m_0+1)^{-1}|\na_X f|^2_{H^1}\le
(\Lambda_a\na_X f,\na_X f)+C(m_0+2)|\na_X f|^2_{L^2}\\
&&\quad \le \f {C(m_0+1)^{-1}} 2|\na_X f|^2_{H^1}+C(m_0+2)|\Lambda_af|^2_{L^2}+C(m_0+2)|f|^2_{L^2}.
\end{eqnarray*}
That is,
\ben\label{2.9}
|\Lambda_af|_{L^2}\ge
C(m_0+2)^{-1}|f|_{H^{2}}-C(m_0+2)|f|_{L^2},
\een
from which and Kato-Ponce commutator estimate, it follows that
\beno
|\Lambda_af|_{H^s}&\ge& |\Lambda_a\Lambda^sf|_{L^2}-|[\Lambda^s,g_{ij}]D_iD_jf|_{L^2}\nonumber\\
&\ge& C(m_0+2)^{-1}|f|_{H^{2+s}}-C(m_0+2)|f|_{H^s}\nonumber\\
&&-C(m_0+2)(|\na_Xa|_{H^{2+s}}|f|_{H^{m_0}}+|f|_{H^2}).
\eeno
Again, by (\ref{2.9}) we get
\ben\label{2.10a}
|\Lambda_af|_{H^s}
\ge C(m_0+2)^{-1}|f|_{H^{2+s}}-C(m_0+2)(1+|\na_Xa|_{H^{2+s}})|f|_{H^{m_0}}.
\een
Now let us inductively assume that for $1\le l\le k-1$
\ben\label{2.10}
|\Lambda_a^lf|_{H^{s}}\ge
C_l(m_0+2)^{-1}|f|_{H^{2l+s}}-C_l(m_0+2)(1+|\na_X a|_{H^{2l+s}})|f|_{H^{m_0}}.
\een
Then we get by the induction assumption that
\ben\label{2.11}
&&|\Lambda_a^kf|_{H^{s}}=|\Lambda_a^{k-1}\Lambda_af|_{H^{s}}
\ge C_{k}(m_0+2)^{-1}|\Lambda_af|_{H^{2(k-1)+s}}\nonumber\\&&\qquad\qquad\qquad
-C_{k}(m_0+2)(1+|\na_X a|_{H^{2(k-1)+s}})|\Lambda_af|_{H^{m_0}}).
\een
While  by (\ref{2.10a}) and Kato-Ponce commutator estimate, we have that for any $|\alpha|=2(k-1)$
\beno
|\partial^\alpha\Lambda_af|_{H^{s}}&\ge&
|\Lambda_a\partial^\alpha f|_{H^{s}}-|[\partial^\alpha,g_{ij}]D_iD_jf|_{H^s}\\
&\ge& C(m_0+2)^{-1}|f|_{H^{2k+s}}-C(m_0+2)|f|_{H^{2(k-1)}}\\
&&-C_k(m_0+2)(|\na_X a|_{H^{2k+s}}|f|_{H^{m_0}}+|f|_{H^{2k-1+s}}),
\eeno
and by interpolation,
\beno
C_k(m_0+2)|f|_{H^{2k-1+s}}\le \f12 C(m_0+2)^{-1}|f|_{H^{2k+s}}+\widetilde{C}_k(m_0+2)|f|_{H^{m_0}},
\eeno
which together with (\ref{2.10}) and (\ref{2.11}) imply the lemma.\ef\vspace{0.2cm}

Similarly, we can also obtain
\bthm{Lemma} \label{lem2.12} Let $k\in \N, s\in [0,1], f\in H^{2k+s}\cap H^{m_0}(\R^d)$. Then there holds
\begin{eqnarray*}
|\sigma_a\Lambda^k_a\,\sigma_a^{-1}f|_{H^{s}}\ge
C_k(m_0+2)^{-1}|f|_{H^{2k+s}}-C_k(m_0+2)
(1+|\na_X
a|_{H^{2k+s}})|f|_{H^{m_0}}.
\end{eqnarray*}
A similar estimate also holds for $|\sigma_a^{-1}\Lambda^k_a\sigma_af|_{H^s}$.
Here $\sigma_a=(1+|\na_X a|^2)^{-\frac14}$.
\ethm

In the energy estimates, we need to deal with the following commutator
\beno
[\cP_a, \Lambda_a^k]=\cP_a\Lambda_a^k-\Lambda_a^k\cP_a,
\eeno
where the operator $\cP_a$ is defined by
\[
\cP_a=(\sigma_a\na_X\,\sigma_a^{-1}\cdot)^2-\Bigl(\frac{\na_Xa\cdot\na_X\sigma_a^{-1}\cdot}{(1+|\na_X
a|^2)^{\frac{3}{4}}}\Bigr)^2.
\]
By a simple calculation, we find that the operator $\cP_a$ can be written as
\begin{eqnarray}\label{2.12}
\cP_af=-\Lambda_af+h_1(\na_Xa,\na^2_Xa)\na_Xf
+h_2(\na_Xa,\na^2_Xa,\na^3_Xa)f,
\end{eqnarray}
for some two smooth functions $h_1,h_2$. So, we arrive at

\bthm{Lemma}\label{lem2.13}
Let $k\in \N, s\in [0,1], f\in H^{2k+s}\cap H^{m_0}(\R^d)$. Then there holds
\[
|[\cP_a,\Lambda^k_a]f|_{H^s}\le
C_k(m_0+4)(|f|_{H^{2k+s}}+|\na_Xa|_{H^{2(k+1)+s}}|f|_{H^{m_0}}).
\]
\ethm

\noindent{\bf Proof.}\, By (\ref{2.12}), we have
\begin{eqnarray*}
[\cP_a,\Lambda^k_a]f=[h_1\na_X,\Lambda^k_a]f+[h_2,\Lambda^k_a]f.
\end{eqnarray*}
Firstly, we get by Lemma \ref{lem2.10} that
\ben\label{2.13}
|[h_2,\Lambda^k_a]f|_{H^s}&\le& |h_2\Lambda^k_af|_{H^s}+|\Lambda^k_a(h_2f)|_{H^s}\nonumber\\
&\le& C_k(m_0+4)(|f|_{H^{2k+s}}+|\na_Xa|_{H^{2(k+1)+s}}|f|_{H^{m_0}}).
\een
While,
\begin{eqnarray*}
&&[h_1\na_X,\Lambda^k_a]f=(h_1\na_X\Lambda^k_af-h_1\Lambda^k_a\na_Xf)
+(h_1\Lambda^k_a\na_Xf-\Lambda^k_ah_1\na_Xf)\\
&&=h_1\sum_{{\sum^k_{l=1}|\alpha_l|+|\beta|=2k+1}\atop{|\beta|\le2k}}
C_{\alpha,\beta}\prod^k_{l=1}D^{\alpha_l}g_{i_lj_l}(\na_Xa)D^\beta f\\
&&\quad+\sum_{{\sum^k_{l=1}|\alpha_l|+|\beta|=2k}}C_{\alpha,\beta}'\prod^k_{l=1}D^{\alpha_l}g_{i_lj_l}(\na_Xa)
\bigl[\sum_{{\beta_1+\beta_2=\beta}\atop{|\beta_2|<|\beta|}}C_{\beta_1,\beta_2}D^{\beta_1}h_1D^{\beta_2}\na_Xf\bigr],\\
\end{eqnarray*}
from which, we get by using H\"{o}lder inequality and interpolation argument that
\ben\label{2.14}
|[h_1\na_X,\Lambda^k_a]f|_{H^s} \le C_k(m_0+2)(|f|_{H^{2k+s}}+|\na_Xa|_{H^{2(k+1)+s}}|f|_{H^{m_0}}).
\een

Summing up (\ref{2.13}) and (\ref{2.14}), we conclude the lemma.\ef

\section{The Dirichlet-Neumann operator}
Assume that the fluid domain $\Omega$  is of the form
$$\Omega=\{(X,y)\in \R^{d+1},b(X)<y<a(X)\}$$
where $a(X)$ and $b(X)$ satisfy (\ref{2.1}).
We consider the boundary value problem
\beq\label{4.1}
\left\{ \begin{array}{l}
-\Delta u=0 \qquad\hbox{on}\quad
\Omega\\
u|_{y=a(X)}=f,\qquad \p_{n_-}u|_{y=b(X)}=0.
\end{array}\right.
\eeq

\bthm{Definition}
Let $k\in\N$, and $a,b\in W^{2,\infty}(\R^d)$ satisfy the condition
(\ref{2.1}). We define the Dirichlet-Neumann operator to be the
operator $G(a,b)$ given by
\[
G(a,b):\begin{matrix} H^{k+\frac32}(\R^d) &
\rightarrow & H^{k+\frac12}(\R^d)\\
f & \mapsto & \sqrt{1+|\na_Xa|^2}\p_{n_+}u|_{y=a(X)},\end{matrix}
\] where $u$ is the solution of (\ref{4.1}).
\ethm

As in section 2, we can associate the elliptic problem (\ref{4.1}) on $\Omega$ to a
problem on a strip $\cS$: \beq\label{4.2} \left\{
\begin{array}{l}
\tilde{\textbf{P}} \tilde{u}=0 \qquad\hbox{on}\quad
\cS,\\
\tilde{u}|_{\tilde{y}=0}=f,\,\,\,\,\,\p^{\tilde{P}}_n\tilde{u}
|_{\tilde{y}=-1}=0.
\end{array}\right.
\eeq We denote by $f^b$ the solution of (\ref{4.2}). Then we have
\[
G(a,b)f=-\p^{\tilde{P}}_nf^b|_{\tilde{y}=0},\,\,\,\,\,\forall f\in
H^\frac32(\R^d).
\]

In what follows, we firstly recall some properties of  the Dirichlet-Neumann operator from \cite{lan1}.
Let us introduce some notations. When a bottom parameterization $b\in W^{k,\infty}(\R^d)$ is given,
we write $B=|b|_{W^{k,\infty}}$. For all $r,s\in \R$, we denote by $M(s)$(resp. $M_r(s))$ constants
which depend on $B$ and $|a|_{H^s}$(resp. $r, B$ and $|a|_{H^s}$).

\bthm{Proposition}\label{Prop3.2}
Assume that $a$, $b$ satisfy (\ref{2.1}). Then there hold

\noindent {\bf i.} For all $k\in \N$, if $a,b \in
W^{k+2,\infty}(\R^d)$, then for all $f$ such that $\na_Xf\in
H^{k+\frac12}(\R^d)^2$, one has
\[
|G(a,b)f|_{H^{k+\frac12}}\le
C(|a|_{W^{k+2,\infty}},|b|_{W^{k+2,\infty}})|\na_Xf|_{H^{k+\frac12}}.
\]

\noindent {\bf ii.} For all $k\in \N$, if $a\in H^{2m_0+\frac12}\cap
H^{k+\frac32}(\R^d)$ and if $b\in W^{k+2}(\R^d)$, then
\[
|G(a,b)f|_{H^{k+\frac12}}\le
M_k(2m_0+\f12)(|\na_Xf|_{H^{k+\frac12}}+|a|_{H^{k+\f32}}
|\na_Xf|_{H^{m_0+\frac12}}),
\]
for all $f$ such that $\na_Xf\in H^{k+\frac12}\cap H^{m_0+\f12}(\R^d)$.
\ethm

\bthm{Proposition}\label{Prop3.3}
Let $a, b\in W^{2,\infty}(\R^d)$ satisfy (\ref{2.1}). Then there hold

\noindent {\bf i.} The operator $G$ is self-adjoint:
\[
(G(a,b)f,g)=(f,G(a,b)g),\quad\forall f,g\in \cS(\R^d);
\]

\noindent {\bf ii.} The operator is positive:
\[
(G(a,b)f,f)\ge 0,\quad \forall f\in \cS(\R^d);
\]

\noindent {\bf iii.} We have the estimates for $f,g\in \cS(\R^d)$:
\begin{eqnarray*}
&&|(G(a,b)f,g)|\le M(m_0+\f12)|f|_{H^\frac12}|g|_{H^\frac12};\\
&&|([G(a,b)+\mu]f,f)|\ge C\tilde{p}|f|_{H^\frac12},
\end{eqnarray*}
for all $\mu\ge\frac{2\tilde{p}}{3}$, where $\tilde{p}$ is given in Lemma \ref{flat-ellip}.
\ethm

Define the operator
\[
R_a=G(a,b)-g_a(X,D),
\]
where $g_a(X,D)$ is a pseudo-differential operator  with the symbol
\[
g_a(X,\xi)=\sqrt{(1+|\na_X a|^2)|\xi|^2-(\xi\cdot \na_Xa)^2}.\]

We have the following tame estimate for $R_a$:

\bthm{Proposition}\label{Prop3.4}
Let $k\in \N, f\in H^{k+\f12}\cap H^{m_0+1}(\R^d)$. Then there holds
\[
|R_af|_{H^{k+\f12}}\le
M_k(2m_0+2)(|f|_{H^{k+\f12}}+|f|_{H^{m_0+1}}|a|_{H^{k+3}}).
\]
\ethm
\bthm{Remark} The case when $k=0,-1$ was proved in \cite{lan1}. To obtain the tame estimate for general $k$,
we need to  use the approximate solution of (\ref{4.2}) constructed in \cite{lan1} and
the tame pseudo-differential operator estimates from \cite{lan2}. Proposition \ref{Prop3.4}  also tells us that
$g_a(X,D)$ is the main part of the Dirichlet-Neumann operator $G(a,b)$.
\ethm

Let us firstly recall the approximate solution of (\ref{4.2}) constructed in \cite{lan1}.
Write $\widetilde{P}$ given by Lemma \ref{flat-ellip} as
\[
\tilde{P}=\left(\begin{matrix} \tilde{P}_1 & \tilde{{\bf p}}\\
\tilde{{\bf p}}^T & \tilde{p}_{d+1}\end{matrix}\right),
\]
i.e.,
\beno \tilde{{\bf
P}}&=&-\tilde{p}_{d+1}\p^2_{\tilde{y}}-(2\tilde{{\bf p}}\cdot
\na_X+(\p_{\tilde{y}}\tilde{p}_{d+1}+\na_X\cdot \tilde{{\bf
p}}))\p_{\tilde{y}}\\
&&-P_1\Delta_X-((\na_X\cdot P_1)+\p_{\tilde{y}}\tilde{{\bf
p}})\cdot\na_X. \eeno
Let
\begin{eqnarray*}
\tilde{\textbf{P}}_{app}=-\tilde{p}_{d+1}(\p_{\tilde{y}}-\eta_-(X,\tilde{y},D))
(\p_{\tilde{y}}-\eta_+(X,\tilde{y},D))
\end{eqnarray*}
where $\eta_\pm(X,\tilde{y},D)$ are pseudo-differential operators with symbols
\[
\eta_\pm(X,\tilde{y},\xi)=\frac 1{\tilde{p}_{d+1}}
\Bigl(-i\tilde{{\bf p}}\cdot\xi\pm\sqrt{\tilde{p}_{d+1}\xi\cdot
\tilde{P}_1\xi-(\tilde{{\bf p}}\cdot\xi)^2}\Bigr),\quad \widetilde{y}\in [-1,0].
\]
Moreover, $\eta_+$ satisfies
\beno
\f {\|\widetilde{P}\|_\infty} {\widetilde{p}}|\xi|\ge \textrm{Re}(\eta_+(X,\widetilde{y},\xi))\ge C_+|\xi|,
\eeno
where $C_+$ is a positive constant depending on $h_0,p,|b|_{1,\infty}$ and $|\na_Xa|_{H^{m_0}}$.
Therefore, there exist functions $\Sigma_\pm(v,\xi)\in C^\infty(\R^{(d+1)^2},\cM^1)$ 1-regular at the origin such that
$\eta_\pm(X,\tilde{y},\xi)=\Sigma_\pm(\widetilde{P}(X,\widetilde{y}),\xi).$

The approximate solution $f^b_{app}$  of (\ref{4.2}) is defined as follows
\ben\label{4.3}
f^b_{app}(X,\widetilde{y})=\sigma_{app}(X,\widetilde{y},D)f,
\een
where $\sigma_{app}(X,\widetilde{y},\xi)=\exp(-\int_{\widetilde{y}}^0\eta_+(X,y',\xi)dy')$.
We know from \cite{lan1} that
\ben\label{4.4}
g_a(X,D)f=-\p^{\tilde{P}}_n
f^b_{app}|_{\tilde{y}=0}.
\een

\noindent{\bf Proof of Proposition \ref{Prop3.4}.}\, By the definition of $G(a,b)$ and (\ref{4.4}), we have
\[
R_af=-\p^{\tilde{P}}_n(f^b-f^b_{app})|_{\tilde{y}=0}\triangleq
-\p^{\tilde{P}}_nf_r^b|_{\tilde{y}=0},
\]
where $f_r^b=f^b-f^b_{app}$. We use the trace theorem and  (\ref{coeff}) to get
\begin{eqnarray*}
|R_af|_{H^{k+\f12}}&=&|e_{d+1}\cdot\tilde{P}\na_{X,\tilde{y}}f_r^b|_{\tilde{y}=0}|_{H^{k+\f12}}
\le C_k\|\tilde{P}\na_{X,\tilde{y}}f_r^b\|_{H^{k,2}}\\
&\le& M_k(m_0+1)(\|f_r^b\|_{k+2,2}+|a|_{H^{k+2}}\|f_r^b\|_{m_0+1,2}).
\end{eqnarray*}

We next turn to the tame estimate of
$\|f_r^b\|_{k+2,2}$. From the proof of Lemma 3.14 in \cite{lan1}, we find that
\[
\tilde{\textbf{P}}f_r^b=-(\tilde{\textbf{P}}-\tilde{\textbf{P}}_{app})
f^b_{app}-\tilde{\textbf{P}}_{app}f^b_{app}:=h^1_{app}+h^2_{app},
\]
together with the boundary condition
$$
f_r^b|_{\widetilde{y}=0}=0,\qquad
\partial_n^{\widetilde{P}}f_r^b|_{\widetilde{y}=-1}=-\partial_n^{\widetilde{P}}f_{app}^b|_{\widetilde{y}=-1}.
$$
From Prop. \ref{tame-elliptic}, it follows that
\begin{eqnarray*}
&&\|f_r^b\|_{k+2,2}\le
M_k(2m_0+2)(\|h^1_{app}\|_{k,2}+\|h^2_{app}\|_{k,2}
+|\p^{\tilde{P}}_nf^b_{app}|_{\tilde{y}=-1}|_{H^{k+\f12}})\\
&&\quad+M_k(2m_0+2)(\|h^1_{app}\|_{m_0-1,2}+\|h^2_{app}\|_{m_0-1,2}
+|\p^{\tilde{P}}_nf^b_{app}|_{\tilde{y}=-1}|_{H^{m_0-\frac12}})|a|_{H^{k+2}}.
\end{eqnarray*}
Now it remains to estimate the right hand side of the above inequality.\vspace{0.1cm}

{\bf Step 1.} Estimate of $h^1_{app}$.\vspace{0.1cm}

We set
\begin{eqnarray*}
&&\tau_1(X,\tilde{y},D)=\eta_-(X,\tilde{y},D)\circ\eta_+(X,\tilde{y},D)
-(\eta_-\eta_+)(X,\tilde{y},D),\\
&&\tau_2(X,\tilde{y},D)=(\p_{\tilde{y}}\eta_+)(X,\tilde{y},D).
\end{eqnarray*}
Then we find that
\begin{eqnarray*}
&&\tilde{\textbf{P}}-\tilde{\textbf{P}}_{app}-\tilde{p}_{d+1}\tau_1
+\tilde{p}_{d+1}\tau_2\\
&&\quad=(\p_{\tilde{y}}\tilde{p}_{d+1}+\na_X\cdot\tilde{{\bf p}})\p_{\tilde{y}}
-(\na_X\cdot\tilde{P}_1+\p_{\tilde{y}}\tilde{{\bf p}})\cdot\na_X.
\end{eqnarray*}
Therefore, we obtain
\begin{eqnarray*}
&&\|(\tilde{\textbf{P}}-\tilde{\textbf{P}}_{app})u\|_{k,2}
\le\|[(\p_{\tilde{y}}\widetilde{p}_{d+1}+\na_X\cdot\tilde{{\bf p}})
\p_{\tilde{y}}-(\na_X\cdot\tilde{P}_1+\p_{\tilde{y}}\tilde{{\bf p}})
\cdot\na_X]u\|_{k,2}\\
&&\qquad\quad+\|\tilde{p}_{d+1}\tau_1(X,\tilde{y},D)u\|_{k,2}
+\|\tilde{p}_{d+1}\tau_2(X,\tilde{y},D)u\|_{k,2}.
\end{eqnarray*}
By (\ref{coeff}), the first term of the right hand side is bounded by
\beno
M_k(m_0+2)(\|u\|_{k+1,2}+|a|_{H^{k+2}}\|u\|_{m_0+1,2}).
\eeno
By Prop. \ref{prop2.8}-\ref{prop2.9} and (\ref{coeff}), the other two terms are bounded by
\beno
M_k(m_0+2)(\|u\|_{k+1,2}+|a|_{H^{k+2}}\|u\|_{m_0+1,2}).
\eeno
So we get
\[
\|(\tilde{\textbf{P}}-\tilde{\textbf{P}}_{app})u\|_{k,2}\le
 M_k(m_0+2)(\|u\|_{k+1,2}+|a|_{H^{k+2}}\|u\|_{m_0+1,2}),
\]
which implies that
\ben\label{4.6}
&&\|h^1_{app}\|_{k,2}=\|(\tilde{\textbf{P}}-\tilde{\textbf{P}}_{app})\widetilde{{\sigma}} _{app}(X,\widetilde{y},D)\exp(\frac{C_+}2 \tilde{y}|D|)f\|_{k,2}\nonumber\\
&&\le M_k(m_0+2)\|\widetilde{{\sigma}} _{app}(X,\widetilde{y},D)\exp(\frac{C_+}2
\tilde{y}|D|)f\|_{k+1,2}\nonumber\\
&&\quad+M_k(m_0+2)\|\widetilde{{\sigma}} _{app}(X,\widetilde{y},D)\exp(\frac{C_+}2
\tilde{y}|D|)f\|_{m_0+1,2}|a|_{H^{k+2}},
\een
where $\widetilde{{\sigma}} _{app}(X,\widetilde{y},D)={\sigma} _{app}(X,\widetilde{y},D)\exp(-\frac{C_+}2
\tilde{y}|D|)$ is a pseudo-differential operator of order 0. By Prop. \ref{prop2.8}, we get
\[
\|\widetilde{{\sigma}} _{app}u\|_{k+1,2}\le
M_k(m_0+1)(\|u\|_{k+1,2}+\|u\|_{m_0,2}|a|_{H^{k+2}}),
\]
which together with (\ref{4.6}) gives
\ben\label{4.7}
\|h^1_{app}\|_{k,2}
\le M_k(m_0+2)(|f|_{H^{k+\f12}}+|f|_{H^{m_0+\f12}}|a|_{H^{k+2}}),
\een
where we used the fact that
\ben\label{4.8}
\|\exp(\frac{C_+}2 \tilde{y}|D|)f\|_{k+1,2}\le C_k|f|_{H^{k+\f12}}.
\een

\vspace{0.1cm}

{\bf Step 2.} Estimate of $h^2_{app}$.\vspace{0.1cm} By the definition, we have
\begin{eqnarray*}
\|h^2_{app}\|_{k,2}
=\|\tilde{p}_{d+1}(\p_{\tilde{y}}-\eta_-(X,\tilde{y},D))\tau_3(X,\tilde{y},D)
\exp(\frac{C_+}2\tilde{y}|D|)f\|_{k,2}.
\end{eqnarray*}
where
$$
\tau_3(X,\tilde{y},D)=(\eta_+\tilde{\sigma}_{app})
(X,\tilde{y},D)-\eta_+\circ\tilde{\sigma}_{app}(X,\tilde{y},D).
$$
We get by Prop. \ref{prop2.8} and (\ref{coeff}) that
\begin{eqnarray*}
\|(\p_{\tilde{y}}-\eta_-(X,\tilde{y},D))u\|_{k,2}
 \le M_k(m_0+1)(\|u\|_{k+1,2}+\|u\|_{m_0+1}|a|_{H^{k+1}}),
\end{eqnarray*}
and by Prop. \ref{prop2.9} and (\ref{coeff})
\beno
\|\tau_3(X,\tilde{y},D)u\|_{k+1,2}\le  M_k(m_0+2)(\|u\|_{k+1,2}+\|u\|_{m_0+1}|a|_{H^{k+3}}),
\eeno
from which and (\ref{4.8}), it follows that
\ben\label{4.9}
\|h^2_{app}\|_{k,2}\le M_k(m_0+2)(|f|_{H^{k+\f12}}+|f|_{H^{m_0+\f12}}|a|_{H^{k+3}}).
\een

\vspace{0.1cm}
{\bf Step 3.} Estimate of $\p^{\tilde{P}}_nf^b_{app}|_{\tilde{y}=-1}$. We rewrite  it as
\begin{eqnarray*}
\p^{\tilde{P}}_nf^b_{app}|_{\tilde{y}=-1}
=e_{d+1}\cdot \tilde{P}|_{\tilde{y}=-1} B(X,D)f+e_{d+1}\cdot
\tilde{P}|_{\tilde{y}=-1}\Bigl(\begin{array}{c}\tau_4 (X,-1,D)\\
0\end{array}\Bigr)f,
\end{eqnarray*}
where
\begin{eqnarray*}
&&\tau_4 (X,\tilde{y},D)=\na_X
\sigma_{app}(X,\tilde{y},D)-(i\xi\sigma_{app})(X,\tilde{y},D),\\
&&B(X,D)=[\Bigl(\begin{array}{l}i\xi\\
\eta_+(X,-1,\xi)\end{array}\Bigr) \sigma_{app}(X,-1,\xi)](X,-1,D).
\end{eqnarray*}
Note that the symbol $\sigma_{app}(X,-1,\xi)$ is a smooth symbol, we get by Prop. \ref{prop2.8} and (\ref{coeff}) that
\beno
|e_{d+1}\cdot \tilde{P}|_{\tilde{y}=-1}B(X,D)f|_{H^{k+\f12}}+|\tau_4(X,-1,D)f|_{H^{k+\f12}}\le
M_k(m_0+2)|a|_{H^{k+\f52}}|f|_{L^2},
\eeno
which leads to
\ben\label{4.10}
|\p^{\tilde{P}}_nf^b_{app}|_{\tilde{y}=-1}|_{H^{k+\f12}}\le
M_k(m_0+2)|a|_{H^{k+\f52}}|f|_{L^2}.
\een

Summing up (\ref{4.7}),(\ref{4.9}) and (\ref{4.10}), we obtain
\ben\label{4.11}
\|f_r^b\|_{k+2,2}\le M_k(2m_0+2)(|f|_{H^{k+\f12}}+|f|_{H^{m_0+1}}|a|_{H^{k+3}}).
\een

The proof of Proposition \ref{Prop3.4} is finished. \ef

\vspace{0.2cm}

We next prove a tame commutator estimate which plays a key role in the energy estimate.
It should be pointed out that we don't need this kind of estimates in the case of zero surface tension.

\bthm{Theorem}\label{prop3.6}
Let $k\in \N, f\in H^{2k+\f12}\cap H^{m_0+2}(\R^d)$. Then there holds
\[
|[\sigma_aG(a,b)\sigma_a,\Lambda_a^k]f|_{H^\f32}\le
M_k(2m_0+2)(|f|_{H^{2k+\f12}}+|f|_{H^{m_0+2}}
|a|_{H^{2k+4}}),
\]
where $\sigma_a=(1+|\na_X a|^2)^{-\frac14}$.
\ethm

\bthm{Remark} The above result seems surprising, since
$[\sigma_aG(a,b)\sigma_a,\Lambda_a^k]$ should be an operator of order $2k$ .
However, thanks to special form of the operator $\Lambda_a$,
the main part of $\sigma_aG(a,b)\sigma_a$ is the same as that of $\Lambda_a$ so that it becomes an operator of order $2k-1$,
which is a key point of this paper.
\ethm
The proof of Theorem \ref{prop3.6} is very technical and will be divided into two parts.
In the first part, we deal with the commutator estimate between the main part of $G(a,b)$ and $\Lambda_a$ which can
be obtained by using pseudo-differential operators calculus. In the second part, we deal with
the commutator estimate between the remainder of $G(a,b)$ and $\Lambda_a$ which relies on the construction
of the approximate solution of the variable coefficients elliptic equation on a flat strip.

\bthm{Lemma}\label{lem3.7}
Let $k\in \N, f\in H^{k+\f32}\cap H^{m_0+2}(\R^d)$. Then there holds
\[
|[\sigma_ag_a(X,D)\sigma_a,\Lambda_a]f| _{H^{k+\f12}}\le
M_k(m_0+3)(|f|_{H^{k+\f32}}+|f|_{H^{m_0+2}}|a|_{H^{k+\f92}}).
\]
\ethm

\noindent{\bf Proof.} We write
\begin{eqnarray*}
&&[\sigma_ag_a(X,D)\sigma_a,\Lambda_a]\\
&&=\sigma_a
g_a(X,D)[\sigma_a,\Lambda_a]+[\sigma_a
g_a(X,D),\Lambda_a]\sigma_a\\
&&=\sigma_ag_a(X,D)\circ([\sigma_a,\Lambda_a]
-\textrm{Op}\{\sigma_a,\Lambda_a\}_1)+\sigma_ag_a(X,D)\circ
\textrm{Op}\{\sigma_a,\Lambda_a\}_1\\
&&\quad+([\sigma_ag_a(X,D),\Lambda_a]-\textrm{Op}\{\sigma_a
g_a(X,\xi),\Lambda_a\}_1)\sigma_a+\textrm{Op}\{\sigma_a
g_a(X,\xi),\Lambda_a\}_1\sigma_a\\
&&\triangleq\sigma_ag(X,D)\circ\tau_{1}(X,D)+\sigma_ag_a(X,D)\circ
\textrm{Op}\{\sigma_a,\Lambda_a\}_1+\tau_{2}(X,D)\sigma_a\\
&&\quad+\textrm{Op}\{\sigma_ag_a(X,\xi),\Lambda_a\}_1\sigma_a.
\end{eqnarray*}
Furthermore, we find that
\begin{eqnarray*}
&&\sigma_ag_a(X,D)\circ
\textrm{Op}\{\sigma_a,\Lambda_a\}_1+ \textrm{Op}\{\sigma_ag_a(X,\xi),\Lambda_a\}_1\sigma_a\\
&&=\bigl[\sigma_ag_a(X,D)\circ \textrm{Op}\{\sigma_a,\Lambda_a\}_1-\textrm{Op}(\sigma_a
g_a(X,\xi)\{\sigma_a,\Lambda_a\}_1)\bigr]\\
&&\quad+\bigl[\textrm{Op}\{\sigma_ag_a(X,\xi),\Lambda_a\}_1\sigma_a-\textrm{Op}(\{\sigma_a
g_a(X,\xi),\Lambda_a\}_1\sigma_a)\bigr]\\
&&\quad +\textrm{Op}(\sigma_a
g_a(X,\xi)\{\sigma_a,\Lambda_a\}_1)+\textrm{Op}(\{\sigma_a
g_a(X,\xi),\Lambda_a\}_1\,\sigma_a)\\
&&\triangleq\tau_{3}(X,D)+\tau_{4}(X,D)
+\textrm{Op}(\{\sigma_ag_a(X,\xi)\sigma_a,\Lambda_a\}_1).
\end{eqnarray*}
Note that
\[
\textrm{Op}(\{\sigma_ag_a(X,\xi)\sigma_a,\Lambda_a\}_1)=0.
\]
Then we obtain
\[
[\sigma_ag_a(X,D)\sigma_a,\Lambda_a]=\sigma_ag(X,D)\circ\tau_{1}(X,D)+\tau_{2}(X,D)\sigma_a
+\tau_{3}(X,D)+\tau_{4}(X,D).
\]
With this identity, the lemma can be deduced from Prop. \ref{prop2.8}-\ref{prop2.9}.
\ef

\bthm{Lemma}\label{lem3.9}
Let $k\in \N, f\in H^{k+\f32}\cap H^{m_0+\f32}(\R^d)$. Then there holds
\[
\|[\tilde{\textbf{P}}-\tilde{\textbf{P}}_{app},\Lambda_a]f^b_{app}\|
_{k,2}\le
M_k(m_0+3)(|f|_{H^{k+\f32}}+|f|_{H^{m_0+\f32}}|a|_{H^{k+4}}).
\]
\ethm

\noindent {\bf Proof.} \, Let $\tau_1, \tau_2$ be as in Prop. \ref{Prop3.4}. We write
\begin{eqnarray*}
[\tilde{\textbf{P}}-\tilde{\textbf{P}}_{app}, \Lambda_a]f^b_{app}&=&
[\tilde{\textbf{P}}-\tilde{\textbf{P}}_{app}+\tilde{p}_{d+1}\tau_2(X,\tilde{y},D)
-\tilde{p}_{d+1}\tau_1(X,\tilde{y},D),\Lambda_a]f^b_{app}\\
&&+[\tilde{p}_{d+1}\tau_1(X,\tilde{y},D), \Lambda_a]f^b_{app}
-[\tilde{p}_{d+1}\tau_2(X,\tilde{y},D),\Lambda_a]f^b_{app}\\
&\triangleq& I_1+I_2+I_3.
\end{eqnarray*}
From the proof of Prop. \ref{Prop3.4}, we find that
\[
I_1=[(\p_{\tilde{y}}\tilde{p}_{d+1}+\na_X\cdot\tilde{{\bf
p}})\p_{\tilde{y}}-(\na_X\cdot\tilde{P}_1+\p_{\tilde{y}}\tilde{{\bf
p}})\cdot\na_X,\Lambda_a]f^b_{app}.
\]
So we get by Prop. \ref{prop2.9} and (\ref{coeff}) that
\begin{eqnarray*}
\|I_1\|_{k,2} \le M_k(m_0+3)(|f|_{H^{k+\f32}}+|f|_{H^{m_0+\frac32}}
|a|_{H^{k+4}}),
\end{eqnarray*}
where we used the fact that(see also (\ref{4.7}))
\beno
\|f^b_{app}\|_{k+2,2}
\le M_k(m_0+1)(|f|_{H^{k+\f32}}+|f|_{H^{m_0+\f12}}|a|_{H^{k+3}}).
\eeno
Similarly, we have
\begin{eqnarray*}
\|I_2\|_{k,2}+\|I_3\|_{k,2} \le M_k(m_0+2)(|f|_{H^{k+\f32}}+|f|_{H^{m_0+\frac32}}
|a|_{H^{k+3}}).
\end{eqnarray*}

This completes the proof of Lemma \ref{lem3.9}.\ef

\bthm{Lemma}\label{lem3.10}
Let $k\in \N, f\in H^{k+\f32}\cap H^{m_0+2}(\R^d)$. Then there holds
\[
|[R_a, \Lambda_a]f|_{H^{k+\f12}}\le
M_k(2m_0+2)(|f|_{H^{k+\f32}}+|f|_{H^{m_0+2}}|a|_{H^{k+5}}).
\]
\ethm

\noindent {\bf Proof.}\, Note that
\[
R_af=-\p^{\tilde{P}}_n(f^b-f^b_{app})|_{\tilde{y}=0}\triangleq
-\p^{\tilde{P}}_nf_r^b|_{\tilde{y}=0},
\]
so we have
\begin{eqnarray*}
[R_a, \Lambda_a]f=e_{d+1}\cdot\tilde{P}\na_{X,\tilde{y}}((\Lambda_af)_r^b
-\Lambda_af_r^b)|_{\tilde{y}=0}+e_{d+1}\cdot[\tilde{P}
\na_{X,\tilde{y}},\Lambda_a]f_r^b|_{\tilde{y}=0}.
\end{eqnarray*}

Firstly, we get by Prop. \ref{prop2.9}, (\ref{coeff}) and (\ref{4.11}) that
\ben\label{4.12}
|[\tilde{P}\na_{X,\tilde{y}},\Lambda_a]f_r^b|_{\tilde{y}=0}|_{H^{k+\f12}}\le M_k(2m_0+2)(|f|_{H^{k+\f32}}+|f|_{H^{m_0+1}}
|a|_{H^{k+4}}).
\een
We now estimate  the term $|e_{d+1}\cdot\tilde{P}\na_{X,\tilde{y}}((\Lambda_af)_r^b
-\Lambda_af_r^b)|_{\tilde{y}=0}|_{H^{k+\f12}}$. By the trace theorem, it suffices to
estimate $\|(\Lambda_af)_r^b -\Lambda_af_r^b\|_{k+2,2}$. By the definition,
\begin{eqnarray*}
&&\tilde{\textbf{P}}f_r^b=-(\tilde{\textbf{P}}-\tilde{\textbf{P}}_{app})f^b_{app}-\tilde{\textbf{P}}_{app}f^b_{app},\\
&&\tilde{\textbf{P}}(\Lambda_af)_r^b=-
(\tilde{\textbf{P}}-\tilde{\textbf{P}}_{app})(\Lambda_af)^b_{app}
-\tilde{\textbf{P}}_{app}(\Lambda_af)^b_{app},
\end{eqnarray*}
which lead to
\begin{eqnarray*}
&&\tilde{\textbf{P}}((\Lambda_af)_r^b
-\Lambda_af_r^b)\\
&&= -[\tilde{\textbf{P}}-\tilde{\textbf{P}}_{app},\Lambda_a]f^b_{app}
+(\tilde{\textbf{P}}-\tilde{\textbf{P}}_{app})(\Lambda_af^b_{app}
-(\Lambda_af)^b_{app})-[\tilde{\textbf{P}},
\Lambda_a]f_r^b\\&&\quad-\tilde{\textbf{P}}_{app}(\Lambda_af)^b_{app}+\Lambda_a\tilde{\textbf{P}}_{app}f^b_{app}
\triangleq h,
\end{eqnarray*}
together with the following boundary conditions
\beno
&&(\Lambda_af)_r^b-\Lambda_af_r^b|_{\tilde{y}=0}=0,\quad \textrm{and} \\
&&-\p^{\tilde{P}}_n((\Lambda_af)_r^b-\Lambda_af_r^b)|
_{\tilde{y}=-1}\\
&&=\bigl(e_{d+1}\cdot\tilde{P}\na_{X,\tilde{y}}(\Lambda_af)^b_{app}-e_{d+1}\cdot\Lambda_a\tilde{P}\na_{X,\tilde{y}}
f^b_{app}+e_{d+1}\cdot[\tilde{P}\na_{X,\tilde{y}},
\Lambda_a]f_r^b\bigr)|_{\tilde{y}=-1}.
\eeno
We get by Prop. \ref{tame-elliptic}  that
\beno
&&\|(\Lambda_af)_r^b
-\Lambda_af_r^b\|_{k+2,2}\\
&&\le M_k(2m_0+2)(\|h\|_{k,2}
+|\p^{\tilde{P}}_n((\Lambda_af)_r^b
-\Lambda_af_r^b)|_{\tilde{y}=-1}|_{H^{k+\f12}})\\
&&\quad+M_k(2m_0+2)(\|h\|_{m_0-1,2}
+|\p^{\tilde{P}}_n((\Lambda_af)_r^b
-\Lambda_af_r^b)|
_{\tilde{y}=-1}|_{H^{m_0-\frac12}})|a|_{H^{k+2}}.
\eeno

Following the proof of Step 3 in Prop. \ref{Prop3.4} and (\ref{4.12}), we find that
\beno
|\p^{\tilde{P}}_n((\Lambda_af)_r^b
-\Lambda_af_r^b)|_{\tilde{y}=-1}|_{H^{k+\f12}}\le M_k(2m_0+2)(|f|_{H^{k+\f32}}+|f|_{H^{m_0+1}}
|a|_{H^{k+4}}).
\eeno

By the definition of $h$, we have
\begin{eqnarray*}
\|h\|_{k,2}&\le& \|[\tilde{\textbf{P}}-\tilde{\textbf{P}}_{app},\Lambda_a]f^b_{app}\|_{k,2}+\|(\tilde{\textbf{P}}-\tilde{\textbf{P}}_{app})(\Lambda_af^b_{app}
-(\Lambda_af)^b_{app})\|_{k,2}\\
&&+\|[\tilde{\textbf{P}},
\Lambda_a]f_r^b\|_{k,2} +\|\Lambda_a
\tilde{\textbf{P}}_{app}f^b_{app}-\tilde{\textbf{P}}_{app}(\Lambda_af)^b_{app}\|
_{k,2}\\
&\triangleq& I+II+III+IV.
\end{eqnarray*}
We get by Lemma \ref{lem3.9} that
\ben\label{4.13}
I\le M_k(m_0+3)(|f|_{H^{k+\f32}}+|f|_{H^{m_0+\f32}}|a|_{H^{k+4}}).
\een
and by the proof of Proposition \ref{Prop3.4},
\ben
&&III \le M_k(2m_0+2)(|f|_{H^{k+\f32}}+|f|_{H^{m_0+1}}
|a|_{H^{k+3}}).\label{4.14}
\een
and
\ben
II&\le& M_k(m_0+2)\bigl(\|\Lambda_af^b_{app}
-(\Lambda_af)^b_{app}\|_{k+1,2}\nonumber\\
&&+\|\Lambda_2f^b_{app}-(\Lambda_af)^b_{app}\|
_{m_0+1,2}|a|_{H^{k+2}}\bigr).\nonumber
\een
Note that
\begin{eqnarray*}
&&\|\Lambda_af^b_{app}
-(\Lambda_af)^b_{app}]\|_{k+1,2}\\
&&\le \|[\Lambda_a,\tilde{\sigma}_{app}]
\exp(\frac{C_+}2\tilde{y}|D|)f\|_{k+1,2}+\|\tilde{\sigma}_{app}
[\Lambda_a,\exp(\frac{C_+}2\tilde{y}|D|)]f\|_{k+1,2}.
\end{eqnarray*}
By Prop. \ref{prop2.9}, (\ref{coeff}) and (\ref{4.8}), the first term  above is bounded by
$$
M_k(m_0+2)(|f|_{H^{k+\f32}}+|f|_{H^{m_0+\f32}}
|a|_{H^{k+4}}).
$$
By the definition of $\widetilde{{\sigma}} _{app}$ in Prop \ref{Prop3.4}, we rewrite it as
$$
\widetilde{{\sigma}} _{app}(X,\widetilde{y},D)={\sigma} _{app}(X,\widetilde{y},D)\exp(-\frac{3C_+}4
\tilde{y}|D|)\exp(\frac{C_+}4
\tilde{y}|D|).
$$
Note that ${\sigma} _{app}(X,\widetilde{y},D)\exp(-\frac{3C_+}4
\tilde{y}|D|)$ is still a pseudo-differential operator of order $0$,  we thus get
 by Prop. \ref{prop2.8}, (\ref{coeff}) and (\ref{4.8}) that
 \ben\label{4.15}
\|\widetilde{{\sigma}} _{app}u(\cdot,\widetilde{y})\|_{k+1}\le M_k(m_0+1)\sup_{\widetilde{y}\in [-1,0]}
(|u(\cdot,\widetilde{y})|_{H^{k+\f12}}+|u(\cdot,\widetilde{y})|_{H^{m_0}}
|a|_{H^{k+2}}).
 \een
On the other hand, we can show by using Bony's paraproduct decomposition\cite{Bon} that for any $\widetilde{y}\in [-1,0]$
\beno
|[\Lambda_a,\exp(\frac{C_+}2\tilde{y}|D|)]f|_{H^{k+\f12}}\le M_k(m_0+2)(|f|_{H^{k+\f32}}+|f|_{H^{m_0+2}}
|a|_{H^{k+4}}),
\eeno
which together with (\ref{4.15}) gives
\ben\label{4.16}
\|\tilde{\sigma}_{app}
[\Lambda_a,\exp(\frac{C_+}2\tilde{y}|D|)]f\|_{k+1,2}\le M_k(m_0+2)(|f|_{H^{k+\f32}}+|f|_{H^{m_0+2}}
|a|_{H^{k+4}}).
\een
So, we get
\ben\label{4.17}
II\le M_k(m_0+2)(|f|_{H^{k+\f32}}+|f|_{H^{m_0+2}}
|a|_{H^{k+4}}).
\een
For $IV$, we have
\[
IV=\|[\Lambda_a,\tilde{\textbf{P}}_{app}\circ\sigma_{app}]f\|_{k,2},
\] where
\begin{eqnarray*}
&&\tilde{\textbf{P}}\circ\sigma_{app}=-\tilde{p}_{d+1}(\p_{\tilde{y}}-\eta_-(X, \tilde{y},D))
\tau(X,\tilde{y},D)\exp(\frac{C_+}2\tilde{y}|D|),\\
&&\tau(X,\tilde{y},D)=(\eta_+\tilde{\sigma}_{app})
(X,\tilde{y},D)-\eta_+\circ\tilde{\sigma}_{app}(X,\tilde{y},D).
\end{eqnarray*}
So, we get
\begin{eqnarray*}
&&[\Lambda_a,\tilde{\textbf{P}}_{app}\circ\sigma_{app}]
=\tilde{p}_{d+1}(\p_{\tilde{y}}-\eta_-)
[\tau,\Lambda_a]\exp(\frac{C_+}2\tilde{y}|D|)\\&&\quad+
[\tilde{p}_{d+1}(\p_{\tilde{y}}-\eta_-),\Lambda_a]
\tau\exp(\frac{C_+}2\tilde{y}|D|)+\tilde{p}_{d+1}(\p_{\tilde{y}}-\eta_-)
\tau[\exp(\frac{C_+}2\tilde{y}|D|),\Lambda_a],
\end{eqnarray*}
from which, we get by  Prop. \ref{prop2.8}-\ref{prop2.9},  (\ref{coeff}) and (\ref{4.8}) (for the last term,
we need to use a similar argument leading to (\ref{4.16})) that
\ben\label{4.18}
IV\le M_k(m_0+3)(|f|_{H^{k+\f32}}+|f|_{H^{m_0+2}}
|a|_{H^{k+5}}).
\een

Summing up (\ref{4.13}-\ref{4.14}) and (\ref{4.17}-\ref{4.18}), we obtain
\beno
\|h\|_{k,2}\le M_k(2m_0+2)(|f|_{H^{k+\f32}}+|f|_{H^{m_0+2}}
|a|_{H^{k+5}}).
\eeno
Thus, we have
\beno
\|(\Lambda_af)_r^b
-\Lambda_af_r^b\|_{k+2,2}
\le M_k(2m_0+2)(|f|_{H^{k+\f32}}+|f|_{H^{m_0+2}}
|a|_{H^{k+5}}),
\eeno
from which and (\ref{4.12}), we conclude the lemma. \ef\vspace{0.2cm}

Now we are in a position to prove Theorem \ref{prop3.6}.\vspace{0.2cm}

\noindent{\bf Proof of Theorem \ref{prop3.6}.}\,\, By the definition of $G(a,b)$,
\begin{eqnarray*}
&&[\sigma_aG(a,b)\sigma_a,\Lambda^k_a]f
=\sum^{k-1}_{i=0}\Lambda^i_a[\sigma_ag_a(X,D)\sigma_a,\Lambda_a]\Lambda^{k-i-1}_af\\
&&\qquad\qquad+\sum^{k-1}_{i=0}\Lambda^i_a[\sigma_a
R_a\sigma_a,\Lambda_a]\Lambda^{k-i-1}_af
\triangleq I+II.
\end{eqnarray*}
Using Lemma \ref{lem2.10} and Lemma \ref{lem3.7}, we get by the interpolation argument that
\beno
|I|_{H^\f32}\le M_k(m_0+3)(|f|_{H^{2k+\f12}}+|f|_{H^{m_0+2}}
|a|_{H^{2k+\f72}}).
\eeno
And by Lemma \ref{lem2.10}, Prop. \ref{Prop3.4}, Prop. \ref{prop2.9} and Lemma \ref{lem3.10}, we have
\beno
|II|_{H^\f32}\le M_k(2m_0+2)(|f|_{H^{2k+\f12}}+|f|_{H^{m_0+2}}
|a|_{H^{2k+4}}).
\eeno

This completes the proof of Theorem \ref{prop3.6}.\ef

\section{The linearized water wave equations}

\subsection{The linearized system}

We rewrite the water wave equations (\ref{1.12}) as
\beq\label{5.1}
\p_tU+\cF(U)=0,
\eeq
where $U=(\zeta,\psi)^T$ and $\cF(U)=(\cF_1(U),\cF_2(U)^T$ with

\begin{equation} \label{5.2}
\begin{split}
\cF_1(U)=&-G(\zeta)\psi, \\
\cF_2(U)=&g\zeta+\frac{1}{2}|\na_X \psi|^2-\frac{(G(\zeta)\psi+\na_X\zeta
\cdot \na_X \psi)^2}{2(1+|\na_X \zeta|^2)}\\
&-\kappa\na_X \cdot (\frac{\na_X
\zeta}{\sqrt{1+|\na_X \zeta|^2}}).
\end{split}
\end{equation}

We will linearize the system (\ref{5.1}) around an admissible reference state in the following sense:

\bthm{Definition}
{Let $T > 0$. We say that
$\underline{U}=(\underline{\zeta},\underline{\psi})^T$ is an
admissible reference state if $(\underline{\zeta},
\underline{\psi}-\underline{\psi}|_{t=0})^T \in C([0,T];
H^\infty(\R^d)^2)$ and $\na_X \underline{\psi}|_{t=0} \in
H^\infty(\R^d)^d$, and if moreover
\beno \exists h_0>0 \qquad
\hbox{such\,\,that}\quad \min\{-b, \underline{\zeta}-b\}\ge h_0 \qquad
\hbox{on}\quad [0,T]\times \R^d,
\eeno
where  $y=b(X)$ is a parameterization of the bottom.}
\ethm

The linearized operator $\underline{\cL}$ associated to
(\ref{5.1}) is given by
 \[ \underline{\cL}:=\p_t+d_{\underline{U}} \cF,
\]
where
\beno
d_{\underline{U}}\cF=\left(
\begin{matrix}
-d_{\underline{\zeta}}G(\cdot)\underline{\psi} & -G(\underline{\zeta})\\
g-\underline{Z}d_{\underline{\zeta}}G(\cdot)\underline{\psi}-\underline{Zv}\cdot
\na_X-\underline{A} & -\underline{Z}G(\underline{\zeta})
+\underline{v}\cdot \na_X
\end{matrix} \right),\eeno
with $\underline{Z}=Z(\underline{U})$,
$\underline{A}=A(\underline{U})$ and
$\underline{v}=v(\underline{U})$, and for all $U=(\zeta,\psi)^T$ smooth enough
\beno
&&Z(U):=\frac{1}{1+|\na_X
\zeta|^2}(G(\zeta)\psi+\na_X
\zeta\cdot \na_X\psi),\label{5.3}\\
&&A(U):=\kappa\na_X\cdot\Bigl[ \frac{\na_X}{\sqrt{1+|\na_X
\zeta|^2}}-\frac{\na_X
\zeta(\na_X\zeta\cdot \na_X)}{(1+|\na_X
\zeta|^2)^\frac{3}{2}}\Bigr], \label{5.4}\\
&&v(U):=\na_X\psi-Z(U)\na_X\zeta.\label{5.5}
\eeno

According to Theorem 3.20 in \cite{lan1}, we have
\[
d_{\underline{\zeta}}G(\cdot)\underline{\psi}\cdot
\zeta=-G(\underline{\zeta})(\underline{Z}\zeta)-\na_X\cdot(\zeta\underline{v}),\]
so that $\underline{\cL}$ becomes
\[ \underline{\cL}=\p_t+\left(
\begin{matrix}
G(\underline{\zeta})(\underline{Z}\cdot)+\na_X\cdot(\cdot\underline{v}) & -G(\underline{\zeta})\\
\underline{Z}G(\underline{\zeta})(\underline{Z}\cdot)+(g+\underline{Z}\na_X\cdot
\underline{v})-\underline{A} &
-\underline{Z}G(\underline{\zeta})+\underline{v}\cdot\na_X
\end{matrix} \right).
\]

Taking $V=(\zeta,\psi-\underline{Z}\zeta)^T$ as a new unknown, the linearized equation
$\underline{\cL}U=G$ is equivalent to
\beq
\underline{\cM}V=H\quad \textrm{on}\quad [0,T]\times\R^d,\label{5.6}
\eeq
where
\[H:=\left(\begin{matrix} G_1\\ G_2-\underline{Z}G_1\end{matrix}\right)
\qquad \hbox{and}\quad\
\underline{\cM}:=\p_t+\left(\begin{matrix}\na_X\cdot(\cdot\underline{v})
& -G(\underline{\zeta})\\
\underline{a}-\underline{A} & \underline{v}\cdot\na_X
\end{matrix}\right) \]
with $\underline{a}:=g+\p_t\underline{Z}+
\underline{v}\cdot\na_X\underline{Z}$.

\subsection{Well-posedness of the linearized system}

We consider the linearized system
\beq \label{5.7} \ \left\{
\begin{array}{ll}
\underline{\cM}V=H \\
V|_{t=0}=V_0
\end{array}\right., \quad \textrm{with}\quad \underline{\cM}=\p_t+\left(\begin{matrix}\na_X\cdot(\cdot\underline{v})
& -G(\underline{\zeta})\\
\underline{a}-\underline{A} & \underline{v}\cdot\na_X
\end{matrix}\right).
\eeq
We obtain the following well-posedness and tame energy estimates of (\ref{5.7}).

\bthm{Proposition} \label{prop5.2} Let $T>0$ and
$\underline{U}$ be an admissible reference state. Assume that $H\in
C([0,T]\times H^\infty(\R^d)^2)$ and $V_0\in H^\infty(\R^d)^2$. Then
there is a unique solution $V\in C^1([0,T],H^\infty(\R^d)^2)$ to
(\ref{5.7}) and for all $k\in \N$, there exists a constant $C_k$
such that
\begin{eqnarray*}
&&|V(t)|_{H^{2k+1}\times H^{2k+\frac{1}{2}}} \\
&&\le C_ke^{C_k t} \bigl[|V_0|_{H^{2k+1}\times H^{2k+\frac{1}{2}}}
+|V_0|_{H^{2m_0+1}\times H^{2m_0+\f12}}(1+|
\underline{\zeta}|_{H^{2k+2}})\bigr] \\
&&+C_k\int^t_0 e^{C_k(t-\tau)}\bigl[ |H|_{H^{2k+1}\times
H^{2k+\frac{1}{2}}}+|H|_{H^{2m_0+1}\times H^{2m_0+\f12}}(1+|
\underline{\zeta}|_{H^{2k+2}})\bigr]d\tau \\
&&+C_k|V_0|_{H^{2m_0+1}\times H^{2m_0+\f12}}\\
&&\quad\times\int^t_0 e^{C_k(t-\tau)}\bigl(1+|\underline{\zeta}|_{H^{2k+4}}+|
\underline{v}|_{H^{2k+2}}+|\underline{a}-g|_{H^{2k+\f12}}+|\p_t
\underline{\zeta}|_{H^{2k+2}}\bigr)d\tau ,
\end{eqnarray*} where
$C_k=C(k,\kappa, B,|\underline{\zeta}|_{H^{q_0}},|
\underline{v}|_{H^{q_0}_T},|\underline{a}-g|_{H^{q_0}_T},|\p_t
\underline{\zeta}|_{H^{q_0}})$, and $q_0$ is a constant depending
only on $d$. Moreover, if $\underline{\cM}$ satisfies the L\'{e}vy condition
\ben\label{levy}
\exists\, c_0>0\quad \textrm{such} \,\,\textrm{that} \quad \underline{a}(t,X)\ge c_0,
\quad \forall (t,X)\in [0,T]\times \R^d,
\een
there exists a constant $\widetilde{C}_k$ independent of $\kappa$
such that
\begin{eqnarray*}
&&|\sqrt{\kappa}V_1(t)|_{H^{2k+1}}+|V(t)|_{H^{2k}\times H^{2k+\frac{1}{2}}}\le C_ke^{C_k t} \bigl[|\sqrt{\kappa}V_{01}|_{H^{2k+1}}
\\&&\quad+|V_{0}|_{H^{2k}\times H^{2k+\frac{1}{2}}}+|V_0|_{H^{2m_0}\times H^{2m_0+\f12}}(1+|
\underline{\zeta}|_{H^{2k+2}})\bigr] \\
&&+C_k\int^t_0 e^{C_k(t-\tau)}\bigl[ |H|_{H^{2k+1}\times
H^{2k+\frac{1}{2}}}+|H|_{H^{2m_0}\times H^{2m_0+\f12}}(1+|
\underline{\zeta}|_{H^{2k+2}})\bigr]d\tau \\
&&+C_k|V_0|_{H^{2m_0}\times H^{2m_0+\f12}}\\
&&\quad\times\int^t_0 e^{C_k(t-\tau)}\bigl(1+|\underline{\zeta}|_{H^{2k+4}}+|
\underline{v}|_{H^{2k+2}}+|\underline{a}-g|_{H^{2k+\f12}}+|\p_t
\underline{\zeta}|_{H^{2k+2}}\bigr)d\tau ,
\end{eqnarray*} where
$\widetilde{C}_k=C(k,c_0, B,|\underline{\zeta}|_{H^{q_0}},|
\underline{v}|_{H^{q_0}_T},|\underline{a}-g|_{H^{q_0}_T},|\p_t
\underline{\zeta}|_{H^{q_0}})$.

\ethm

\noindent{\bf Proof.}\, As the existence of solutions to (\ref{5.7})
follows from the {\it a priori} estimates for the approximate
solutions( which can be constructed by a parabolic regularization as
in \cite{lan1}), here we only  present the {\it a priori} tame
estimate to smooth enough solutions of (\ref{5.7}).

We rewrite the linear system (\ref{5.7}) as
\beq\label{5.8}  \left\{
\begin{array}{ll}
\p_tV_1+\na_X\cdot(\underline{v}V_1)-\underline{G}V_2=H_1 \\
\p_tV_2+(\underline{a}-\underline{A})V_1+\underline{v}\cdot \na_X V_2=H_2\\
V|_{t=0}=V_0
\end{array}\right.
\eeq where we write $\underline{G}$ for $G(\underline{\zeta})$. We introduce the following energy functional $E_k(V)$ defined by
\[ \noindent
E_k(V):=(\widetilde{\Lambda}^k\sigma V_1
,\sigma^{-1}(\underline{a}-\underline{A})\,\sigma^{-1}\widetilde{\Lambda}^k\sigma
V_1)+(\widetilde{\Lambda}^k\sigma^{-1}V_2,\sigma \underline{G}
\,\sigma\widetilde{\Lambda}^k \sigma^{-1} V_2),
\]
where
\beno
&&\widetilde{\Lambda}=|D|^2-\frac{\p_i\underline{\zeta}\p_j\underline{\zeta}}
{1+|\na_X \underline{\zeta}|^2}D_iD_j,\quad\sigma=(1+|\na_X
\underline{\zeta}|^2)^{-\frac{1}{4}}.
\eeno

Fix a constant $\lambda>0$ to be determined later. We have
\begin{eqnarray}\label{5.9}
\noindent &&\frac{d}{dt}e^{-2\lambda t}E_k(V)\nonumber\\
&&= e^{-2\lambda t}\bigl[-2\lambda
E_k(V)+2(\tilde{\Lambda}^k\,\p_t(\,\sigma
V_1),\sigma^{-1}(\underline{a}-\underline{A})\sigma^{-1}\tilde{\Lambda}^k\,\sigma
V_1)\nonumber\\
&&\quad+2(
\tilde{\Lambda}^k\p_t(\sigma^{-1}V_2),\sigma \underline{G}
\sigma \tilde{\Lambda}^k\sigma^{-1} V_2)+( \tilde{\Lambda}^k\sigma
V_1,[\p_t,\,\sigma^{-1}(\underline{a}-\underline{A})\sigma^{-1}]\tilde{\Lambda}^k\sigma
V_1)\nonumber\\
&&\quad+( \tilde{\Lambda}^k\sigma^{-1}
V_2,[\p_t,\,\sigma \underline{G}\sigma]\tilde{\Lambda}^k
\sigma^{-1}V_2)+2( [\p_t,\,\tilde{\Lambda}^k]\sigma
V_1,\sigma^{-1}(\underline{a}-\underline{A})\sigma^{-1}\tilde{\Lambda}^k\sigma
V_1)\nonumber\\
&&\quad +2([\p_t,\tilde{\Lambda}^k]\,\sigma^{-1}
V_2,\sigma \underline{G} \,\sigma \tilde{\Lambda}^k
\,\sigma^{-1} V_2)\bigr]\nonumber\\
&&\triangleq e^{-2\lambda t}(I_1+\cdots+I_7).\nonumber
\end{eqnarray}
Let us begin with the estimates of $I_1-I_7$.\vspace{0.2cm}

{\bf Estimtates of $I_2+I_3$.}
\begin{eqnarray*}
I_2& =&2(\tilde{\Lambda}^k(\p_t\,\sigma)
V_1,\sigma^{-1}(\underline{a}-\underline{A})\sigma^{-1}\tilde{\Lambda}^k\sigma
V_1)\\&&+2(\tilde{\Lambda}^k\sigma\p_t
V_1,\sigma^{-1}(\underline{a}-\underline{A})\sigma^{-1}\tilde{\Lambda}^k\sigma
V_1)\\
&\triangleq& I_{21}+I_{22}.
\end{eqnarray*}
By Lemma \ref{lem2.10}, we  have
\ben\label{5.10}
|I_{21}|\le C_k[|V_1|^2_{H^{2k}}+\kappa |V_1|^2_{H^{2k+1}}+|V_1|^2_{H^{m_0}}
(|\underline{\zeta}|^2_{H^{2k+2}}
+|\p_t\underline{\zeta}|^2_{H^{2k+2}})].
\een
Using the first equation of (\ref{5.8}), we rewrite $I_{22}$ as
\begin{eqnarray*}
I_{22}&=& 2(\tilde{\Lambda}^k\sigma\underline{G}V_2,\sigma^{-1}(\underline{a}-\underline{A})\sigma^{-1}\tilde{\Lambda}^k\sigma
V_1)\\&&-2(\tilde{\Lambda}^k\sigma\na_{\underline{v}}
V_1,\sigma^{-1}(\underline{a}-\underline{A})\sigma^{-1}\tilde{\Lambda}^k\sigma
V_1)\\
&&-2(\tilde{\Lambda}^k\sigma\frac{1}{2}(\na_X\cdot
\underline{v})V_1,\sigma^{-1}(\underline{a}-\underline{A})\sigma^{-1}\tilde{\Lambda}^k\sigma
V_1)\\&&+2(\tilde{\Lambda}^k\sigma
H_1,\sigma^{-1}(\underline{a}-\underline{A})\,\sigma^{-1}\tilde{\Lambda}^k\sigma
V_1)\\
&\triangleq& I_{22}^1+\cdots+I_{22}^4.
\end{eqnarray*}
Here the operator $\na_{\underline{v}}$ is defined by
\[
\na_{
\underline{v}}f\triangleq\frac{1}{2}(\na_X\cdot(f\,\underline{v})+\underline{v}\cdot
\na_X f).
\]
Since $\na_{\underline{v}}$ is an anti-adjoint operator, we have
\beno
I_{22}^2&=& 2([\tilde{\Lambda}^k\sigma,\na_{\underline{v}}]
V_1,\sigma^{-1}(\underline{a}-\underline{A})\sigma^{-1}\tilde{\Lambda}^k\sigma
V_1)\\&&-(\tilde{\Lambda}^k\sigma
V_1,[\sigma^{-1}(\underline{a}-\underline{A})\sigma^{-1},\na_{\underline{v}}]
\tilde{\Lambda}^k\sigma V_1),
\eeno
which together with Lemma \ref{lem2.10} gives
\ben\label{5.11}
|I_{22}^2|\le C_k\bigl[|V_1|^2_{H^{2k}}+\kappa |V_1|^2_{H^{2k+1}}+|V_1|^2_{H^{m_0}}
(|\underline{\zeta}|^2_{H^{2k+3}}
+|\underline{v}|^2_{H^{2k+2}})\bigr].
\een
Again, we get by Lemma \ref{lem2.10} that
\begin{eqnarray}
|I_{22}^3|\le C_k\bigl[|V_1|^2_{H^{2k}}+\kappa |V_1|^2_{H^{2k+1}}+|V_1|^2_{H^{m_0}}
(|\underline{\zeta}|^2_{H^{2k+2}}
+|\underline{v}|^2_{H^{2k+2}})\bigr],\label{5.12}
\end{eqnarray}
and
\ben
|I_{22}^4|&\le& C_k\bigl[|H_1|^2_{H^{2k+1}}+|V_1|^2_{H^{2k}}+\kappa |V_1|^2_{H^{2k+1}}\nonumber\\
&&\qquad+(|H_1|^2_{H^{m_0}}+|V_1|^2_{H^{m_0}})
|\underline{\zeta}|^2_{H^{2k+2}}\bigr].\label{5.13}
\end{eqnarray}
The term $I_{22}^1$ will be handled together with $I_3$. We have
\begin{eqnarray*}
I_3&=&2(\tilde{\Lambda}^k(\p_t\sigma^{-1})V_2,\sigma \underline{G}
\sigma \tilde{\Lambda}^k\sigma^{-1} V_2)+2( \tilde{\Lambda}^k\sigma^{-1}\p_tV_2,\sigma
\underline{G}\sigma \tilde{\Lambda}^k\sigma^{-1}V_2)\\
&\triangleq&I_{31}+I_{32}.
\end{eqnarray*}
We get by Prop. \ref{Prop3.3} and Lemma \ref{lem2.10} that
\begin{eqnarray}
|I_{31}|&\le& C_k(|\sigma\tilde{\Lambda}^k\p_t(\sigma^{-1})V_2|
^2_{H^{\frac12}}+|\sigma \tilde{\Lambda}^k\sigma^{-1} V_2|^2_{H^{\frac12}})\nonumber\\
&\le& C_k\bigl[|V_2|^2_{H^{2k+\frac12}}+|V_2|^2_{H^{m_0}}
(|\underline{\zeta}|^2_{H^{2k+\frac32}}
+|\p_t\underline{\zeta}|^2_{H^{2k+\frac32}})\bigr].\label{5.14}
\end{eqnarray}
From the second equation of (\ref{5.8}), we get
\begin{eqnarray*}
I_{32}&=&-2(
\tilde{\Lambda}^k\sigma^{-1}(\underline{a}-\underline{A})V_1,\sigma
\underline{G}\sigma \tilde{\Lambda}^k\sigma^{-1}
V_2)\\
&&-2(
\tilde{\Lambda}^k\sigma^{-1}\na_{\underline{v}}V_2,\sigma
\underline{G}\sigma \tilde{\Lambda}^k\sigma^{-1}
V_2)\\&&+2(
\tilde{\Lambda}^k\sigma^{-1}\frac{1}{2}(\na_X\cdot
\underline{v})V_2,\sigma \underline{G}\sigma
\tilde{\Lambda}^k\sigma^{-1} V_2)\\
&&+2(
\tilde{\Lambda}^k\sigma^{-1}H_2,\sigma \underline{G}\sigma \tilde{\Lambda}^k\sigma^{-1} V_2)\\
&\triangleq&I_{32}^1+\cdots+I_{32}^4.
\end{eqnarray*}
Again, by the fact that $\na_{\underline{v}}$ is anti-adjoint, we have
\[
I_{32}^2=-2([\tilde{\Lambda}^k\sigma^{-1},\na_{\underline{v}}]
V_2,\sigma \underline{G}\sigma \tilde{\Lambda}^k\sigma^{-1}V_2)
+(\tilde{\Lambda}^k\sigma^{-1}V_2,[\sigma \underline{G}\sigma,
\na_{\underline{v}}]\tilde{\Lambda}^k\sigma^{-1}
V_2).
\]
 The first term above can be handled as in $I_{31}$, and for the
second one we need to use Prop. 3.18 in \cite{lan1}. Then we get by Lemma \ref{lem2.10} that
\ben
|I_{32}^2|\le C_k\bigl[|V_2|^2_{H^{2k+\frac12}}+|V_2|^2_{H^{m_0}}
(|\underline{\zeta}|^2_{H^{2k+\f52}}
+|\underline{v}|^2_{H^{2k+\f32}})\bigr].\label{5.16}
\een
Similarly, we have
\ben
&&|I_{32}^3|\le C_k\bigl[|V_2|^2_{H^{2k+\frac12}}+|V_2|^2_{H^{m_0}}
(|\underline{\zeta}|^2_{H^{2k+\f32}}
+|\underline{v}|^2_{H^{2k+\f32}})\bigr],\label{5.17}\\
&&|I_{32}^4|\le C_k\bigl[|H_2|^2_{H^{2k+\frac12}}+|V_2|^2_{H^{2k+\frac12}}
+(|H_2|^2_{H^{m_0}}
+|V|^2_{H^{m_0}})|\underline{\zeta}|^2_{H^{2k+\f32}}\bigr].\label{5.18}
\een

Now it remains to estimate $I_{22}^1+I_{32}^1$. We have
\begin{eqnarray*}
I_{22}^1+I_{32}^1&=& 2(\tilde{\Lambda}^k\sigma\underline{G}V_2,\sigma^{-1}(\underline{a}-\underline{A})\sigma^{-1}\tilde{\Lambda}^k\sigma
V_1)\\&&-2(
\tilde{\Lambda}^k\sigma^{-1}(\underline{a}-\underline{A})V_1,\sigma
\underline{G}\sigma \tilde{\Lambda}^k\sigma^{-1}V_2)\\
&\triangleq& 2(II+III).
\end{eqnarray*}
By the definition of $\underline{A}$,
\beno
\underline{A}&=&\kappa\na_X\cdot\Bigl[ \frac{\na_X}{\sqrt{1+|\na_X
\underline{\zeta}|^2}}-\frac{\na_X
\underline{\zeta}(\na_X\underline{\zeta}\cdot \na_X)}{(1+|\na_X
\underline{\zeta}|^2)^\frac{3}{2}}\Bigr]\\
&=&\kappa\na_X\cdot[\sigma^2\na_X-\sigma^6\na_X
\underline{\zeta}(\na_X\underline{\zeta}\cdot \na_X)],
\eeno
so we get
\begin{eqnarray*}
II&=& (\tilde{\Lambda}^k\sigma\underline{G}V_2,\underline{a}\sigma^{-2}\tilde{\Lambda}^k\sigma
V_1)\\
&&-\kappa(\tilde{\Lambda}^k\sigma\underline{G}V_2,\sigma^{-1}\na_X\cdot\
\sigma^2\na_X\sigma^{-1}\tilde{\Lambda}^k\sigma
V_1)\\
&&
+\kappa(\tilde{\Lambda}^k\sigma\underline{G}V_2,\sigma^{-1}\na_X\cdot\sigma^6\na_X
\underline{\zeta}(\na_X\underline{\zeta}\cdot
\na_X\sigma^{-1}\tilde{\Lambda}^k\sigma
V_1))\\
&\triangleq& II_1+II_2+II_3.
\end{eqnarray*}
We rewrite $II_1$ as
\beno
II_1&=&(\sigma\underline{G}\sigma\tilde{\Lambda}^k\sigma^{-1}V_2,\tilde{\Lambda}^k\underline{a}\sigma^{-1}
V_1)+([\tilde{\Lambda}^k,\sigma\underline{G}\sigma]\sigma^{-1}V_2,\underline{a}\sigma^{-2}\tilde{\Lambda}^k\sigma
V_1)\nonumber\\
&&+(\sigma\underline{G}\sigma\tilde{\Lambda}^k\sigma^{-1}V_2,[\underline{a}\sigma^{-2},\tilde{\Lambda}^k]\sigma
V_1)\nonumber\\
&\triangleq& II_{11}+II_{12}+II_{13}.
\eeno
By Lemma \ref{lem2.10}, Prop. \ref{Prop3.2} and Thm. \ref{prop3.6}, we have
\ben\label{5.14a}
|\cR_1|\triangleq |II_{12}+II_{13}|\le C_k\bigl[|V_1|^2_{H^{2k}}+|V_2|^2_{H^{2k+\f12}}+|V|^2_{H^{m_0+2}} |\underline{\zeta}|^2_{H^{2k+4}}\bigr].
\een
Now, we decompose $II_2$ as
\begin{eqnarray*}
II_2&
=&-\kappa(\tilde{\Lambda}^k\sigma\underline{G}V_2,\sigma^{-1}\na_X\cdot
\sigma [\sigma\na_X\sigma^{-1},\,\tilde{\Lambda}^k]\sigma
V_1)\\&&-\kappa(\tilde{\Lambda}^k\,\sigma\underline{G}V_2,\sigma^{-1}\na_X\cdot\,\sigma\tilde{\Lambda}^k\sigma\na_X
V_1)\\
&=&-\kappa(\tilde{\Lambda}^k\sigma\underline{G}V_2,\sigma^{-1}(\na_X\sigma)\cdot[\sigma\na_X\sigma^{-1},
\tilde{\Lambda}^k]\sigma V_1)\\
&&-\kappa(\tilde{\Lambda}^k\sigma\underline{G}V_2,\sigma^{-1}\sigma\na_X\cdot[\sigma\na_X\sigma^{-1},
\tilde{\Lambda}^k]\sigma V_1)\\
&&-\kappa(\tilde{\Lambda}^k\sigma\underline{G}V_2,\sigma^{-1}(\na_X\sigma)\cdot\tilde{\Lambda}^k\sigma\na_X
V_1)\\&&-\kappa(\tilde{\Lambda}^k\sigma\underline{G}V_2,\sigma^{-1}\sigma\na_X\cdot\tilde{\Lambda}^k\sigma\na_X
V_1)\\
&\triangleq&II_{21}+II_{22}+II_{23}+II_{24}.
\end{eqnarray*}
We get by Lemma \ref{lem2.10} and Prop. \ref{Prop3.2} that
\ben
|II_{21}|\le
C_k\bigl[\kappa |V_1|^2_{H^{2k+1}}+|V_2|^2_{H^{2k+\f12}}+|V|^2_{H^{m_0}} |\underline{\zeta}|^2_{H^{2k+3}}\bigr].\label{5.19}
\een
We write $II_{22}$ as
\begin{eqnarray*}
II_{22}&
=&-\kappa(\tilde{\Lambda}^k\sigma\underline{G}V_2,\na_X\cdot[\sigma\na_X\,\sigma^{-1},\tilde{\Lambda}^k]\sigma
V_1)\\
&=&-\kappa(\tilde{\Lambda}^k\sigma\underline{G}V_2,\sigma\na_X\cdot\sigma^{-1}[\sigma\na_X\,\sigma^{-1},\tilde{\Lambda}^k]\sigma
V_1)\\
&&+\kappa(\tilde{\Lambda}^k\sigma\underline{G}V_2,\sigma(\na_X\sigma^{-1})\cdot[\sigma\na_X\,\sigma^{-1},\tilde{\Lambda}^k]\sigma
V_1)\\
&\triangleq& II_{22}^1+II_{22}^2.
\end{eqnarray*}
As in (\ref{5.19}), we have
\ben
|II_{22}^2|\le
C_k\bigl[\kappa |V_1|^2_{H^{2k+1}}+|V_2|^2_{H^{2k+\f12}}+|V|^2_{H^{m_0}} |\underline{\zeta}|^2_{H^{2k+3}}\bigr].\label{5.20}
\een
We write $II_{23}$ as
\begin{eqnarray*}
II_{23}
&=&-\kappa(\tilde{\Lambda}^k\sigma\underline{G}V_2,[\sigma^{-1}(\na_X\,\sigma),\tilde{\Lambda}^k]\cdot\sigma\na_X
V_1)\\&&-\kappa(\tilde{\Lambda}^k\sigma\underline{G}V_2,\tilde{\Lambda}^k\sigma^{-1}(\na_X\,\sigma)\cdot\sigma\na_X
V_1)\\
&\triangleq&II_{23}^1+II_{23}^2.
\end{eqnarray*}
As in (\ref{5.19}), we have
\ben
|II_{23}^1|\le
C_k\bigl[\kappa |V_1|^2_{H^{2k+1}}+|V_2|^2_{H^{2k}}+|V|^2_{H^{m_0}} |\underline{\zeta}|^2_{H^{2k+3}}\bigr].\label{5.21}
\een
We write $II_{24}$ as
\begin{eqnarray*}
II_{24}&=&-\kappa(\tilde{\Lambda}^k\sigma\underline{G}V_2,\na_X\cdot\tilde{\Lambda}^k\sigma\na_X
V_1)\\
&=&-\kappa(\tilde{\Lambda}^k\sigma\underline{G}V_2,\sigma\na_X\cdot\sigma^{-1}\tilde{\Lambda}^k\sigma\na_X
V_1)\\&&+\kappa(\tilde{\Lambda}^k\sigma\underline{G}V_2,\sigma(\na_X\sigma^{-1})\cdot\tilde{\Lambda}^k\sigma\na_X
V_1)\\
&
=&-\kappa(\tilde{\Lambda}^k\,\sigma\underline{G}V_2,\,\,[\sigma\na_X\,\sigma^{-1},\,\tilde{\Lambda}^k]\cdot\,\sigma\na_X
V_1)\\&&-\kappa(\tilde{\Lambda}^k\,\sigma\underline{G}V_2,\,\,\tilde{\Lambda}^k\sigma\na_X\cdot\na_X
V_1)\\
&&+\kappa(\tilde{\Lambda}^k\sigma\underline{G}V_2,[\sigma(\na_X\sigma^{-1}),\tilde{\Lambda}^k]\cdot\sigma\na_X
V_1)\\&&+\kappa(\tilde{\Lambda}^k\sigma\underline{G}V_2,\tilde{\Lambda}^k\sigma(\na_X\sigma^{-1})\cdot\sigma\na_X
V_1)\\
&\triangleq& II_{24}^1+\cdots+II_{24}^4.
\end{eqnarray*}
As in (\ref{5.19}), we have
\ben
|II_{24}^3|\le
C_k\bigl[\kappa |V_1|^2_{H^{2k+1}}+|V_2|^2_{H^{2k}}+|V|^2_{H^{m_0}} |\underline{\zeta}|^2_{H^{2k+3}}\bigr].\label{5.22}
\een
And from the above calculations, we find that
\begin{eqnarray*}
&&II_{23}^2+II_{24}^2+II_{24}^4\\
&&=-\kappa(\tilde{\Lambda}^k\sigma\underline{G}V_2,\tilde{\Lambda}^k[\sigma^{-1}(\na_X\sigma)\cdot\sigma\na_X
V_1+\sigma\na_X\cdot\na_X
V_1-\sigma(\na_X\sigma^{-1})\cdot\sigma\na_X V_1])\\
&&
=-\kappa(\tilde{\Lambda}^k\sigma\underline{G}V_2,\tilde{\Lambda}^k[\sigma^{-1}(\na_X\sigma)\cdot\sigma\na_X
V_1+\sigma^{-1}\sigma\na_X\cdot\sigma\na_X V_1])\\
&&=-\kappa(\tilde{\Lambda}^k\sigma\underline{G}V_2,\tilde{\Lambda}^k\sigma^{-1}\na_X\cdot\sigma^2\na_X
V_1),
\end{eqnarray*}
and we also have
\begin{eqnarray*}
II_{22}^1+II_{24}^1
&=&-\kappa(\tilde{\Lambda}^k\sigma\underline{G}V_2,\sigma\na_X\cdot\sigma^{-1}[\sigma\na_X\sigma^{-1},\tilde{\Lambda}^k]\sigma
V_1)\\&&-\kappa(\tilde{\Lambda}^k\sigma\underline{G}V_2,[\sigma\na_X\sigma^{-1},\tilde{\Lambda}^k]\cdot\sigma\na_X
V_1)\\
&=&-\kappa(\tilde{\Lambda}^k\,\sigma\underline{G}V_2,\,\,[(\sigma\na_X\,\sigma^{-1}\cdot)^2,\,\tilde{\Lambda}^k]\,\sigma
V_1).
\end{eqnarray*}

Summing up (\ref{5.19})-(\ref{5.22}), we obtain
\begin{eqnarray*}
II_2
&=&-\kappa(\tilde{\Lambda}^k\sigma\underline{G}V_2,\tilde{\Lambda}^k\sigma^{-1}\na_X\cdot(\sigma^2\na_XV_1))\\
&&-\kappa(\tilde{\Lambda}^k\sigma\underline{G}V_2,[(\sigma\na_X\sigma^{-1}\cdot)^2,\tilde{\Lambda}^k]\sigma
V_1)+ {\cR}_2,
\end{eqnarray*}
where the remainder terms $\cR_2$ satisfies
\ben
\cR_2\le
C_k\bigl[|(\sqrt{\kappa} V_1,V_2)|^2_{H^{2k+1}\times{H^{2k+\f12}}}+|V|^2_{H^{m_0}}|\underline{\zeta}|^2_{H^{2k+3}}\bigr].\label{5.23}
\een

Similarly, $II_3$ can be written as
\begin{eqnarray*}
II_3&=&\kappa(\tilde{\Lambda}^k\sigma\underline{G}V_2,\tilde{\Lambda}^k\sigma^{-1}\na_X\cdot\sigma^6\na_X\underline{\zeta}
(\na_X\underline{\zeta}\cdot\na_X V_1))\\
&&
\kappa(\tilde{\Lambda}^k\sigma\underline{G}V_2,[(\sigma^3\na_X\underline{\zeta}\cdot\na_X\sigma^{-1}\cdot)^2,\tilde{\Lambda}^k]\,\sigma
V_1)+\cR_3,
\end{eqnarray*}
where $\cR_3$ has the same estimate as $\cR_2$. So, we arrive at
\begin{eqnarray}
II&
=&(\sigma\underline{G}\sigma\tilde{\Lambda}^k\sigma^{-1}V_2,\tilde{\Lambda}^k\underline{a}\sigma^{-1}
V_1)-(\tilde{\Lambda}^k\sigma\underline{G}V_2,\tilde{\Lambda}^k\sigma^{-1}\underline{A}\,V_1)\nonumber
\\
&&-\kappa(\tilde{\Lambda}^k\sigma\underline{G}V_2,
[(\sigma\na_X\sigma^{-1}\cdot)^2-(\sigma^3\na_X\underline{\zeta}\cdot\na_X\sigma^{-1}\cdot)^2,
\tilde{\Lambda}^k]\,\sigma V_1)\nonumber\\
&&+\cR_1+\cR_2+\cR_3.\label{5.20a}
\end{eqnarray}
On the other hand,
\beno
III&=&-(\tilde{\Lambda}^k\sigma^{-1}\underline{a}V_1,\sigma
\underline{G}\sigma\tilde{\Lambda}^k\sigma^{-1} V_2)+(\tilde{\Lambda}^k\sigma^{-1}\,\underline{A}V_1,[\sigma
\underline{G}\sigma,\tilde{\Lambda}^k]\sigma^{-1} V_2)\\&&+(
\tilde{\Lambda}^k\sigma^{-1}\underline{A}V_1,\tilde{\Lambda}^k\sigma
\underline{G}\,V_2),
\eeno
which together with  (\ref{5.20a}) gives
\begin{eqnarray*}
I_{22}^1+I_{32}^1&=&-2\kappa(\tilde{\Lambda}^k\sigma\underline{G}V_2,
[(\sigma\na_X\sigma^{-1}\cdot)^2-(\sigma^3\na_X\underline{\zeta}\cdot\na_X\sigma^{-1}\cdot)^2,\tilde{\Lambda}^k]
\sigma V_1)\\
&&+2(\tilde{\Lambda}^k\sigma^{-1}\underline{A}V_1,[\sigma
\underline{G}\sigma,\tilde{\Lambda}^k]\sigma^{-1}
V_2)\\
&&+\cR_1+\cR_2+\cR_3\\
&\triangleq& \cR_4+\cR_5+\cR_1+\cR_2+\cR_3.
\end{eqnarray*}
By Lemma \ref{lem2.10}, Lemma \ref{lem2.13} and Thm. \ref{prop3.6},  we have
\begin{eqnarray*}
\cR_4&\le& 2\kappa|\tilde{\Lambda}^k\,\sigma\underline{G}V_2|_{H^{-1}}
|[(\sigma\na_X\sigma^{-1}\cdot)^2-(\sigma^3\na_X\underline{\zeta}\cdot\na_X\sigma^{-1}\cdot)^2,\tilde{\Lambda}^k]\sigma
V_1|_{H^1}\\
&\le&C_k\bigl[|(\sqrt{\kappa} V_1,V_2)|^2_{H^{2k+1}\times
H^{2k+\frac12}}+|V|^2_{H^{m_0}\times H^{m_0}}
|\underline{\zeta}|^2_{H^{2k+4}}\bigr],\\
\cR_5&\le&2|\tilde{\Lambda}^k\sigma^{-1}\underline{A}V_1|_{H^{-1}}
|[\sigma
\underline{G} \sigma,\tilde{\Lambda}^k]\sigma^{-1} V_2|_{H^\f32}\\
&\le& C_k\bigl[|(\sqrt{\kappa} V_1,V_2)|^2_{H^{2k+1}\times
H^{2k+\frac12}}+|V|^2_{H^{m_0}\times H^{m_0+2}}
|\underline{\zeta}|^2_{H^{2k+4}}\bigr].
\end{eqnarray*}
which together with (\ref{5.10})-(\ref{5.18}), (\ref{5.14a}) and and (\ref{5.23}) give
\begin{eqnarray*}
&&|I_2+I_3|\le C_k(\kappa |V_1|^2_{H^{2k+1}}+|V|^2_{H^{2k}\times H^{2k+\frac12}})\\
&&\,\,+C_k|V|^2_{H^{m_0}\times H^{m_0+2}}
\bigl(|\underline{\zeta}|^2_{H^{2k+4}}+|\underline{v}|
^2_{H^{2k+2}}+|\underline{a}-g|_{H^{2k+\f12}}^2+|\p_t\underline{\zeta}|^2_{H^{2k+2}}\bigr)\\
&&\,\,+C_k\bigl[|H|^2_{H^{2k+1}\times
H^{2k+\frac12}}+|H|^2_{H^{m_0}}
|\underline{\zeta}|^2_{H^{2k+2}}\bigr].
\end{eqnarray*}

{\bf Estimates of $I_4-I_7$.} Since the estimates are very
similar as above but much simpler (by Lemma \ref{lem2.10}, and Prop. 3.19 in \cite{lan1} for $I_4$),
we omit it here. We have
\beno
&&|I_4+\cdots+I_7|\\&&\quad\le C_k\bigl[\kappa |V_1|^2_{H^{2k+1}}+|V|^2_{H^{2k}\times H^{2k+\frac12}}+|V|^2_{H^{m_0}}
\bigl(|\underline{\zeta}|^2_{H^{2k+4}}+|\p_t\underline{\zeta}|^2_{H^{2k+2}}\bigr)\bigr].
\eeno

{\bf The total energy estimate.} We finally obtain
\begin{eqnarray}\label{5.24}
&&\frac{d}{dt}e^{-2\lambda t}E_k(V)\nonumber\\
&&\le-2\lambda e^{-2\lambda t}E_k(V)+e^{-2\lambda
t}C_k\bigr(\kappa |V_1|^2_{H^{2k+1}}+|V|^2_{H^{2k}\times H^{2k+\frac12}}\bigr)\nonumber\\
&&+e^{-2\lambda t}C_k|V|^2_{H^{m_0}\times H^{m_0+2}}
(|\underline{\zeta}|^2_{H^{2k+4}}+|\underline{v}|
^2_{H^{2k+2}}+|\underline{a}-g|_{H^{2k+\f12}}^2+|\p_t\underline{\zeta}|^2_{H^{2k+2}}\bigr)\nonumber\\
&&+ e^{-2\lambda t}C_k\bigl[|H|^2_{H^{2k+1}\times
H^{2k+\frac12}}+|H|^2_{H^{m_0}}
|\underline{\zeta}|^2_{H^{2k+2}}\bigr].
\end{eqnarray}

To complete the energy estimates, we still need the following lemma:
\bthm{Lemma}\label{lem5.3} There exists a positive constant $C_{m_0}$ depending only $k$,
$|\underline{a}|_{L^\infty}$, and $|\underline{\zeta}|_{H^{m_0+3}}$  such that the following inequalities hold
\beno
E_k(V)&\ge& C_{m_0}^{-1}(|(\sqrt{\kappa}V_1,V_2)|^2_{H^{2k+1}\times H^{2k+\frac12}}+\inf_{(t,X)\in [0,T]\times \R^d} \underline{a}(t,X)|V_1|^2_{H^{2k}})\\
&&-C_{m_0}|V|^2_{H^{m_0}}(1+|\na_X\underline{\zeta}|^2_{H^{2k+1}}),\\
E_k(V)&\le& C_{m_0}(|(\sqrt{\kappa}V_1,V_2)|^2_{H^{2k+1}\times H^{2k+\frac12}}+|V_1|^2_{H^{2k}})\\&&+C_{m_0}|V|^2_{H^{m_0}}
(1+|\na_X\underline{\zeta}|^2_{H^{2k+1}}).
\eeno

\ethm

\noindent {\bf Proof of Lemma \ref{lem5.3}.} We  write
\[
E_k(V)=
E_k(V_1)+E_k(V_2).
\]
where
\beno
E_k(V_1)=(\tilde{\Lambda}^k\sigma V_1
,\sigma^{-1}(\underline{a}-\underline{A})\sigma^{-1}\tilde{\Lambda}^k\sigma
V_1),\quad
E_k(V_2)=(\tilde{\Lambda}^k\sigma^{-1}V_2,\sigma \underline{G}\sigma \tilde{\Lambda}^k \sigma^{-1}V_2).
\eeno
Set $f=\sigma^{-1}\tilde{\Lambda}^k\sigma
V_1$, we get
\begin{eqnarray*}
E_k(V_1)&=&(f,\underline{a}-\underline{A}f)\\
&=& \int_{\R^d}\underline{a}|f|^2dX+\kappa\int_{\R^d}\frac{(1+|\na_X\underline{\zeta}|^2)|\na_Xf|^2
-|\na_X\underline{\zeta}\cdot\na_Xf|^2}{(1+|\na_X\underline{\zeta}|^2)
^{\frac32}}dX,
\end{eqnarray*}
which implies that
\[
C^{-1}_{m_0}\kappa|\na_Xf|^2_{L^2}\le E_k(V_1)-\int_{\R^d}\underline{a}|f|^2dX\le C_{m_0}\kappa|\na_Xf|^2_{L^2}.
\]
On the other hand, by Lemma \ref{lem2.10}-\ref{2.13} we have
\begin{eqnarray*}
&&|f|^2_{L^2}
\ge C^{-1}_{m_0}|V_1|^2_{H^{2k}}-C_{m_0}
|V_1|^2_{H^{m_0}}(1+|\na_X\underline{\zeta}|^2_{H^{2k}}),\\
&&|f|^2_{L^2}
\le C_{m_0}\bigl[|V_1|^2_{H^{2k}}+
|V_1|^2_{H^{m_0}}(1+|\na_X\underline{\zeta}|^2_{H^{2k}})\bigr],\\
&&|\na_Xf|^2_{L^2}
\ge C^{-1}_{m_0}|V_1|^2_{H^{2k+1}}-C_{m_0}
|V_1|^2_{H^{m_0}}(1+|\na_X\underline{\zeta}|^2_{H^{2k+1}}),\\
&&|\na_Xf|^2_{L^2}
\le C_{m_0}\bigl[|V_1|^2_{H^{2k+1}}+
|V_1|^2_{H^{m_0}}(1+|\na_X\underline{\zeta}|^2_{H^{2k+1}})\bigr],
\end{eqnarray*}
from which, we get
\begin{eqnarray*}
&&E_k(V_1)+C_{m_0}
|V_1|^2_{H^{m_0}}(1+|\na_X\underline{\zeta}|^2_{H^{2k+1}})
\ge C_{m_0}^{-1}\bigl[\kappa|V_1|^2_{H^{2k+1}}+\inf\underline{a}|V_1|^2_{H^{2k}}\bigr],\\
&& E_k(V_1)\le C_{m_0}\bigl[\kappa|V_1|^2_{H^{2k+1}}+|V_1|^2_{H^{2k}}+
|V_1|^2_{H^{m_0}}(1+|\na_X\underline{\zeta}|^2_{H^{2k+1}})\bigr].
\end{eqnarray*}

We now turn to $E_k(V_2)$. Set $g=\sigma\tilde{\Lambda}^k \sigma^{-1}
V_2$, we have
\[
E_k(V_2)=(\tilde{\Lambda}^k\sigma^{-1}V_2,\sigma
\underline{G}\sigma \tilde{\Lambda}^k\sigma^{-1}
V_2)\triangleq(g,\underline{G}g),
\]
which together with Lemma \ref{lem2.10}, Lemma \ref{lem2.12} and
Prop. \ref{Prop3.3} gives
\begin{eqnarray*}
E_k(V_2)&\ge& C|g|^2_{H^{\frac12}}-\mu |g|^2_{L^2}\\
&=& C|\tilde{\Lambda}^k \,\sigma^{-1}V_2|^2_{H^{\frac12}}
-\mu |\tilde{\Lambda}^k \,\sigma^{-1}V_2|^2_{L^2}\\
&\ge& C_{m_0}^{-1}|V_2|^2_{H^{2k+\frac12}}-C_{m_0}
|V_2|^2_{H^{m_0}}(1+|\na_X\underline{\zeta}|^2_{H^{2k+\frac12}})\\
&&-C_{m_0}|V_2|^2_{H^{2k}}-C_{m_0}
|V_2|^2_{H^{m_0}}|\na_X\underline{\zeta}|^2_{H^{2k}}]\\
&\ge& C_{m_0}^{-1}|V_2|^2_{H^{2k+\frac12}}-C_{m_0}
|V_2|^2_{H^{m_0}}(1+|\na_X\underline{\zeta}|^2_{H^{2k+\frac12}}),
\end{eqnarray*} which leads to
\[
E_k(V_2)+C_{m_0}
|V_2|^2_{H^{m_0}}(1+|\na_X\underline{\zeta}|^2_{H^{2k+\frac12}}) \ge
C_{m_0}^{-1}|V_2|^2_{H^{2k+\frac12}}.
\]
And we also have from Prop. \ref{Prop3.3} and Lemma \ref{lem2.10}
\[
E_k(V_2)\le C_{m_0}|V_2|^2_{H^{2k+\frac12}}+ C_{m_0}
|V_2|^2_{H^{m_0}}(1+|\na_X\underline{\zeta}|^2_{H^{2k+\frac12}}).
\]
This finishes the proof of Lemma \ref{lem5.3}.\ef\vspace{0.1cm}

Now we are in a position to complete the proof of Proposition \ref{prop5.2}.
By (\ref{5.24}) and Lemma \ref{lem5.3}, we get
\begin{eqnarray}
&&\frac{d}{dt}e^{-2\lambda t}E_k(V)\nonumber\\
&&\le (-2\lambda+C_k)e^{-2\lambda t}E_k(V)\nonumber\\
&&+e^{-2\lambda t}C_k|V|^2_{H^{m_0}\times H^{m_0+2}}
\bigl(1+|\underline{\zeta}|^2_{H^{2k+4}}+|\underline{v}|
^2_{H^{2k+2}}+|\underline{a}-g|_{H^{2k+\f12}}^2+|\p_t\underline{\zeta}|^2_{H^{2k+2}}\bigr)\nonumber\\
&&+ e^{-2\lambda t}C_k\bigl[|H|^2_{H^{2k+1}\times
H^{2k+\frac12}}+|H|^2_{H^{m_0}}
(1+|\underline{\zeta}|^2_{H^{2k+2}})\bigr].\label{energy}
\end{eqnarray}
Taking $\lambda$ such that $-2\lambda+C_k\le 0$, the first estimate of the Proposition follows easily from  Lemma \ref{lem5.3}.
Note that the constant $C_k$ in (\ref{5.24}) is independent of $\kappa$. So, if $\underline{a}$ satisfies (\ref{levy}),
$C_k$ in (\ref{energy}) is also independent of $\kappa$ by Lemma \ref{lem5.3}.

This completes the proof of Proposition \ref{prop5.2}.\ef

\subsection{The L\'{e}vy condition}

Let us introduce the pressure $\underline{P}$ as
\beno
-\underline{P}=\p_t \underline{\phi}+\f12 |\na_{X,y}\underline{\phi}|^2+gy.
\eeno
where $\underline{\phi}$ is a solution of the elliptic equations
\beno
\left\{\begin{array}{ll}
&-\Delta_{X,y} \underline{\phi}=0 \qquad \textrm{on}\quad \Omega(t),\\
&\underline{\phi}|_{y=\underline{\zeta}(t,X)}=-\underline{\psi}(t,x),
\quad \p_{n_-}\underline{\phi}|_{y=b(X)}=0.
\end{array}\right.
\eeno
If the admissible reference $\underline{U}$ solves the water wave problem (\ref{5.1}) at time $t_0$, we find
from the proof of Proposition 4.4 in \cite{lan1} that
\ben
&&\underline{P}|_{y=\underline{\zeta}(t_0,X)}=-\kappa \na_X\cdot (\frac{\na_X
\underline{\zeta}(t_0,X)}{\sqrt{1+|\na_X \underline{\zeta}(t_0,X)|^2}}),\label{pressure}\\
&&-\p_{n_+}\underline{P}|_{y=\underline{\zeta}(t_0,X)}=\underline{a}(t_0,X)\nonumber\\
&&\qquad-\kappa\frac{\na_X
\underline{\zeta}(t_0,X)\cdot\na_X}{\sqrt{1+|\na_X \underline{\zeta}(t_0,X)|^2}}\na_X\cdot (\frac{\na_X
\underline{\zeta}(t_0,X)}{\sqrt{1+|\na_X \underline{\zeta}(t_0,X)|^2}}).\label{5.24a}
\een
So, the L\'{e}vy condition is equivalent to the Taylor sign condition for $\kappa=0$. In what follows, we show that
\bthm{Proposition}
Let $b=-1$. If $\underline{U}$ solves the water wave problem (\ref{5.1}) at time $t_0$
and the surface tension coefficient $\kappa$ is sufficiently small, the L\'{e}vy condition still holds:
there exists $c_0>0$ such that
\beno
\underline{a}(t_0,X)\ge c_0\quad \textrm{on} \quad \R^d.
\eeno
\ethm

\no{\bf Proof.}\, Note that $-\Delta_{X,y} \underline{P}=\Delta(\f12|\na_{X,y}\underline{\phi}|^2)$ and
\beno
\Delta(\f12|\na_{X,y}\underline{\phi}|^2)=|\na_{X,y}^2\underline{\phi}|^2,
\quad -\p_{n_-}(\f12|\na_{X,y}\underline{\phi}|^2)|_{y=-1}=0,
\eeno
we have
\beno
\left\{\begin{array}{ll}
&-\Delta_{X,y} \underline{P}=|\na_{X,y}^2\underline{\phi}|^2\textrm{on}\quad \Omega(t),\\
&\underline{P}|_{y=\underline{\zeta}(t_0,X)}=-\kappa \na_X\cdot (\frac{\na_X
\underline{\zeta}(t_0,X)}{\sqrt{1+|\na_X \underline{\zeta}(t_0,X)|^2}}),
\quad \p_{n_-}\underline{P}|_{y=-1}=g.
\end{array}\right.
\eeno
From (\ref{5.24a}), it suffices to prove that if $\kappa$ is sufficiently small,
there exists a constant $c_0$ such that $-\p_{n_+}\underline{P}|_{y=\underline{\zeta}(t_0,X)}\ge 2c_0$.
We will follow the argument of Lemma 4.1 in \cite{Wu2}. We denote by $\Gamma_t$ the free surface, and by $\Gamma_b$ the bottom.
For any $h\in C_0^1(\Gamma_t), h\ge 0$, let $w$ be a solution of the elliptic equation
\beno
-\Delta w=0,\quad w|_{\Gamma_t}=h,\quad \f {\p w} {\p n}|_{\Gamma_b}=0.
\eeno
Applying Green's second identity to $w$ and $\underline{P}+gy$, we obtain
\beno
-\int_{\Gamma_t}w\f {\p \underline{P}} {\p n}dS&=&\int_{\Omega(t)}|\na^2_{X,y}\underline{\phi}|^2wdV-\int_{\Gamma_t}\underline{P}\f {\p w}{\p n}dS\\
&&+g\int_{\Gamma_t}(w \f {\p y}{\p n}-y\f {\p w}{\p n})dS.
\eeno
Applying Green's second identity again to $w$ and $y$, we get
\beno
\int_{\Gamma_t}(w \f {\p y}{\p n}-y\f {\p w}{\p n})dS=\int_{\Gamma_b}w dS.
\eeno
Let $G=G(\eta,\xi)$ be the Green's function on $\Omega(t)$, i.e.
\beno
\Delta_\xi G(\eta,\xi)=\delta_\eta(\xi)\quad G(\eta,\xi)|_{\Gamma_t}=0,\quad \f {\p G(\eta,\xi)} {\p n(\xi)}|_{\Gamma_b}=0.
\eeno
Then we have Green's representation formula
\beno
w(\eta)=\int_{\Gamma_t}h(\xi)\f {\p G(\eta,\xi)} {\p n(\xi)}dS(\xi),\quad \textrm{for}\,\, \eta\in \Omega(t).
\eeno
Therefore, we obtain
\ben\label{5.23a}
&&-\int_{\Gamma_t}h(\xi)\f {\p \underline{P}} {\p n}dS(\xi)\ge
\int_{\Gamma_t}h(\xi)\bigl(\int_{\Gamma_b}\f {\p G(\eta,\xi)} {\p n(\xi)}dS(\eta)\bigr)dS(\xi)\nonumber\\
&&\qquad\qquad-\int_{\Gamma_t}h(\xi)\f{\p}{\p n(\xi)}\bigl(\int_{\Gamma_t}\underline{P}(\eta)\f {\p G(\eta,\xi)} {\p n(\eta)}dS(\eta)\bigr)dS(\xi).
\een
From the maximal principle, we know that there exists a constant $c_0$ such that
\beno
\int_{\Gamma_b}\f {\p G(\eta,\xi)} {\p n(\xi)}dS(\eta)\ge 3c_0, \quad \textrm{for any}\,\, \xi\in \Gamma_t.
\eeno
Since $h\ge 0$ is arbitrary, we get from (\ref{5.23a}) and (\ref{pressure}) that
\beno
-\p_{n_+}\underline{P}|_{y=\underline{\zeta}(t_0,X)}\ge 2c_0,
\eeno
if $\kappa$ is sufficiently small. \ef

\section{The Nonlinear Equations and zero surface tension limit}

\subsection{The water wave equations with surface tension}

Assume that the surface tension coefficient $\kappa>0$. Then we prove

\bthm{Theorem}\label{thm6.1}
Let $b\in C^\infty_b(\R^d)$. There exists $P>D>0$ such that for all $\zeta_0\in H^{s+P}(\R^d)$ and $\psi_0$
be such that  $\na_X \psi_0\in H^{s+P}(\R^d)^d$, with $s>M$ ($M$ depends
on $d$ only). Assume moreover that
\[
\min\{\zeta_0-b,-b\}\ge 2h_0 \qquad\hbox{on}\quad \R^d \qquad\hbox{for
some} \quad h_0\ge 0.
\]
Then there exist $T>0$ and a unique solution $(\zeta,\psi)$ to
the water-wave system (\ref{1.12}) with the initial condition
$(\zeta_0,\psi_0)$ and such that $(\zeta,\psi-\psi_0)\in
C^1([0,T],\,H^{s+D}(\R^d)\times H^{s+D}(\R^d))$.
\ethm

We will use a simplified Nash-Moser theorem in \cite{Alv2} to solve (\ref{1.12}). We rewrite the water wave equations (\ref{1.12}) as
\beq\label{6.1}
\p_tU+\cF(U)=0, \quad U(0,x)=U_0(x),
\eeq
where $U=(\zeta,\psi)^T$ and $\cF(U)=(\cF_1(U),\cF_2(U)^T$ with $\cF_1(U)$ and $\cF_2(U)$  given by (\ref{5.2}).
In order to use Theorem 1(with $\varepsilon=1, \cL=0$) in \cite{Alv2},
we firstly verify the assumptions on the nonlinear operator $\cF(\cdot)$. Assume that
$\underline{U}$ is an admissible reference state. Note that $\underline{U}$ dose not necessarily belong to Sobolev space
because $\underline{\psi}|_{t=0}$ is not necessarily in a Sobolev space
(though its gradient is). In what follows, we use $|\underline{U}|_{H^a}$ to
denote the quantity
\[
|\underline{U}|_{H^a}:=|\underline{U}-\underline{U}|_{t=0}|_{H^a}+
|\na_X \underline{U}|_{t=0}|_{H^{a-1}}.
\]
By Prop. \ref{Prop3.2}, there holds for $a\ge q_0$
\begin{eqnarray*}
|\cF(\underline{U})|_{H^{a}_T}\le
C(a,|\underline{U}|_{H_T^{q_0}})|\underline{U}|_{H^{a+2}_T}.
\end{eqnarray*}
We know from the expression of $d_{\underline{U}}\cF$ that
\[
d_{\underline{U}}\cF\cdot U=\left(
\begin{array}{l}
-d_{\underline{\zeta}}G(\cdot)\underline{\psi}\cdot\zeta
-G(\underline{\zeta})\psi\\
g\zeta-\underline{Z}d_{\underline{\zeta}}G(\cdot)
\underline{\psi}\cdot\zeta-\underline{Z}\underline{v}\cdot\na_X\zeta
+\underline{A}\zeta-\underline{Z}G(\underline{\zeta})\psi+\underline{v}
\cdot\na_X\psi
\end{array}\right),
\]
from which and Prop. 3.25 in \cite{lan1}, it follows that for $a\ge q_0$
\[
|d_{\underline{U}}\cF\cdot U| _{H^{a}_T}\le
C(a,|\underline{U}|_{H_T^{q_0}},T)\bigl[|U|_{H^{a+2}_T}
+|U|_{H^{q_0}_T} |\underline{U}|_{H^{a+2}_T}\bigr].
\]
We can compute $d^2_{\underline{U}}\cF$ from the expression of $d_{\underline{U}}\cF$
, and then get by Prop. 3.25 in \cite{lan1} that for $a\ge q_0$
\beno
|d_{\underline{U}}^2\cF\cdot(U,V)| _{H^{a}_T}&\le&
C(a,B,|\underline{U}|_{H_T^{q_0}},T)\bigl[|U|_{H^{q_0}_T} |V|_{H^{a+2}_T}+|U|_{H^{a+2}_T} |V|_{H^{q_0}_T}\\&&
+|U|_{H^{q_0}_T}|V|_{H^{q_0}_T}|\underline{U}|_{H^{a+2}_T}\bigr].
\eeno

To verify the assumptions on the linearized equations,
we need to introduce two scales of Banach spaces $X^a$ and $F^a$ as follows:
\beno
&&X^a=\bigcap_{j=0}^1 C^j([0,T]; H^{a-2j}(\R^d)^2),\quad |u|_{X^a}=\sum_{j=0}^1|\p_t^ju|_{H^{a-2j}_T},\\
&&F^a=C([0,T]; H^{a}(\R^d)^2)\times H^{a}(\R^d)^2,\quad |(f,g)|_{F^a}=|f|_{H_T^a}+|g|_{H^{a}},
\eeno
and for $(f,g)\in F^a$ and $t\in [0,T]$,
\beno
\cI^a(t,f,g)=|g|_{H^{a}}+\int_0^t\sup_{t''\in [0,t']}|f(t'')|_{H^a}dt'.
\eeno
We consider the linearized system of (\ref{6.1})
\beq \label{6.2} \ \left\{
\begin{array}{ll}
\underline{\cM}V=H \\
V|_{t=0}=V_0
\end{array}\right., \quad \textrm{with}\quad \underline{\cM}=\p_t+\left(\begin{matrix}\na_X\cdot(\cdot\underline{v})
& -G(\underline{\zeta})\\
\underline{a}-\underline{A} & \underline{v}\cdot\na_X
\end{matrix}\right).
\eeq
For any $s\ge q_0,\underline{ U}\in X^{s+4}, (H,V_0)\in F^{s+1}$, from Prop. \ref{prop5.2},
there exists a unique solution $V\in C([0,T]; H^s(\R^d))$ to (\ref{6.2}) such that
\ben\label{6.3}
|V|_{H_T^s}\le C(s,\kappa,|\underline{U}|_{X^{q_0}},T)\bigl(\cI^{s+1}(t,H,V_0)+|\underline{U}|_{X^{s+4}}\cI^{q_0}(t,H,V_0)\bigr).
\een
Then Theorem 1 in \cite{Alv2} ensures that there exist $T>0$ and a unique solution $U$
to (\ref{6.1}).

\subsection{Zero surface tension limit} For the flat bottom, we prove that

\bthm{Theorem}\label{thm6.2} Let $b=-1$ and $(\zeta_0, \psi_0)$ satisfy the same assumptions as Theorem \ref{thm6.1}.
Assume  that the surface tension coefficient $\kappa$ is sufficiently small.
Then there exist $T>0$ independent of $\kappa$  and a unique solution $(\zeta^\kappa,\psi^\kappa)$ to
(\ref{1.12}) with $\kappa>0$  such that $(\zeta^\kappa,\psi^\kappa-\psi_0)\in
C^1([0,T],\,H^{s+D}(\R^d)\times H^{s+D}(\R^d))$.
Moreover, as $\kappa$ tends to zero, the solution $(\zeta^\kappa,\psi^\kappa)$ converges to the solution
$(\zeta,\psi)$ of (\ref{1.12}) with $\kappa=0$.
\ethm

\no{\bf Proof.} If the surface tension coefficient $\kappa$ is sufficiently small,
the L\'{e}vy condition holds. Thus, the estimate (\ref{6.3}) is independent of $\kappa$ by Prop \ref{prop5.2}.
By Theorem 1 in \cite{Alv2}, there exist $T>0$ independent of $\kappa$ and a unique solution
$U^\kappa$ of (\ref{6.1}) which is bounded in $C([0,T];H^{s+D}(\R^d))$. Then there exists a subsequence
of $\{U^\kappa\}_{\kappa>0}$ which converges weakly to some $U$. By a standard compactness argument, we can prove that
$U$ is a solution of (\ref{6.1}) with $\kappa=0$.

\section*{Acknowledgements}

{\small We would like  to thank Professor Ping Zhang  and Chongchun Zeng for helpful
discussions. This work was done when Zhifei Zhang was visiting Department of Mathematics of Paris-Sud
University as a Postdoctor Fellowship. He would like to thank the hospitality and support of the Department.}

\end{document}